\documentclass[final, 11pt]{article}

\usepackage{amsfonts}
\usepackage{makeidx}
\usepackage{amsmath}
\usepackage{amssymb}
\usepackage{bm}
\usepackage{cite}
\usepackage{showkeys}
\usepackage[english]{babel}
\usepackage{dblaccnt}
\usepackage{accents}
\usepackage{graphicx}
\usepackage{psfrag}
\usepackage{subfig}
\usepackage{color}
\usepackage{float}
\usepackage{amscd}
\usepackage[mathscr]{euscript}
\usepackage{lipsum}
\usepackage{epstopdf}
\usepackage{algorithm}
\usepackage{algpseudocode}

\newcommand{\eps}{\varepsilon}
\newcommand{\x}{{\bf x}}

\renewcommand{\d}{\mathrm{d}}

\newtheorem{defi}{Definition}

\newtheorem{remark}[defi]{Remark}

\begin{document}

\title{Numerical Solution of a Coefficient Inverse Problem with 
Multi-Frequency Experimental  Raw Data by a Globally
Convergent Algorithm}
\author{Dinh-Liem Nguyen\thanks{
Department of Mathematics and Statistics, University of North Carolina at
Charlotte, Charlotte, NC 28223, USA; (\texttt{dnguye70@uncc.edu}, 
\texttt{mklibanv@uncc.edu}, \texttt{lnguye50@uncc.edu},
 \texttt{akolesov@uncc.edu}, \texttt{hliu34@uncc.edu}) }
\and Michael V. Klibanov\footnotemark[1] \and Loc H. Nguyen\footnotemark[1]  \and 
Aleksandr E. Kolesov\footnotemark[1] \thanks{Institute of Mathematics and Information Science, North-Eastern Federal University, Yakutsk, Russia; 
(\texttt{ae.kolesov@s-vfu.ru})} \and Michael A. Fiddy\thanks{
Optoelectronics Center, University of North Carolina at Charlotte,
Charlotte, NC 28223, USA; (\texttt{mafiddy@uncc.edu})} \and Hui Liu\footnotemark[1] }

\date{}
\maketitle

\begin{abstract}
We analyze in this paper the performance of a newly developed globally convergent numerical method for a
coefficient inverse problem for the case of multi-frequency experimental
backscatter data associated to a single incident wave. These data were collected
using a microwave scattering facility at the University of North Carolina at
Charlotte. The challenges for the inverse problem under the consideration
are not only from its high nonlinearity and severe ill-posedness but also
from the facts that the amount of the measured data is minimal and that
these raw data are contaminated by a significant amount of noise, due to a non-ideal
experimental setup. This setup is motivated by our target application in
detecting and identifying explosives. We show in this paper how the
raw data can be preprocessed and successfully inverted using our
inversion method. More precisely, we are able to reconstruct the dielectric
constants and the locations of the scattering objects with a good accuracy,
without using any advanced
\emph{a priori} knowledge of their physical and geometrical properties.
\end{abstract}


\sloppy


\bigskip




%
%




\section{Introduction}

We are interested in a Coefficient Inverse Problem (CIP) with real world applications, including the detection and identification of explosives, nondestructive testing and material characterization. A new globally convergent
numerical method for solving such a CIP has been recently developed by our group in \cite{Kliba2016}. Numerical study in \cite{Kliba2016} was conducted for
computationally simulated data. In this paper, we study the performance of this method for the case of
experimental raw data. These data were collected using a microwave scattering
facility at the University of North Carolina at Charlotte.


More precisely, we consider the CIP of reconstruction of physical and geometrical properties
of three-dimensional objects from experimental multi-frequency data without
using any detailed \textit{a priori} knowledge of those objects.
Our study is mainly motivated by
potential applications in detection and identification of explosives such
as, e.g., anti-personnel mines and improvised explosive devices (IEDs), see~%
\cite{Schub2006, Weath2015}.
 These targets are
placed in air, and the measured data are the backscatter corresponding to a
single incident wave at multi-frequencies. 
In addition, we note that IEDs are also often buried under
the ground for which one needs to study the corresponding inverse problem of
determining buried objects. This topic is subject to a future publication.

The idea here is to determine the dielectric constants of targets. The
knowledge of dielectric constants might serve in the future classification
algorithms as a piece of information, which would be an additional one to
those commonly currently used in the radar community. Indeed, this community
relies now only on the intensity of radar images, see, e.g.,~\cite{Kuzhu2012, Soumekh:1999}. It is well known that the question of a reliable differentiation between
explosives and clutter is not yet addressed satisfactory in the radar
community. So, estimates of dielectric constants and shapes of targets
combined with the image intensity and other parameters might lead in the
future to such classification algorithms, which would address this question
better. The second potential application of our study is in nondestructive
testing and materials characterization.

A Coefficient Inverse Problem is the problem of recovering a
coefficient of a partial differential equation from boundary measurements of its solutions.
There is a large body of the literature on imaging methods for
reconstructing geometric information about targets, such as their shapes,
sizes and locations. We refer to, e.g.,~\cite{Ammar2004, AmmariChowZou:sjap2016, Burge2005, Cakon2006,
Colto1996, Kirsc1998, LiLiuWang:jcp2014, LiLiuZou:smms2014, Li2015, Pasto2010} and
references therein for well-known imaging techniques in inverse scattering
such as level set methods, sampling methods, expansion methods, and shape optimization methods.
However, for detecting or identifying, for instance, IEDs, the physical properties of the targets (in our model, the
dielectric constants), would play a more important role~\cite{Weath2015}.
Furthermore, determining the spatially distributed dielectric constants,
which is of our main interest, is known to be a more difficult task since
 CIPs are, in general, highly nonlinear and severely ill-posed.

Among the inversion methods developed for solving the CIPs, the two probably most well known approaches are nonlinear
approximation schemes and weak scattering approximation methods such as Born
approximation and physical optics. The weak scattering approximation is not
applicable to the inverse problem under consideration, where the scattering
objects can be strong scatterers. Regarding the nonlinear optimization
approaches or also known as iterative solution methods, we refer to, e.g.,~%
\cite{Engl1996, Chave2010, Gonch2013} and references therein. It is
well-known that the convergence for this class of methods typically requires
a good \textit{a priori} initial approximation of the exact solution, that is, the
starting point of iterations should be chosen to be sufficiently close to
the solution. Hence, we call such methods \textit{locally convergent}. We note that
in our desired applications such \textit{a priori} knowledge is not always available.

This limitation of nonlinear optimization approaches is avoided in the
so-called approximately globally convergent method (globally convergent
method, for short, or GCM), which has been recently introduced, see \cite%
{Beili2012}. The GCM, which does not use optimization schemes, aims to
provide a good approximation to the solution of the coefficient inverse
problem without using any advanced \textit{a priori} knowledge of the solution. More
precisely, the concept ``approximate globally convergence" can be understood
in the language of functional analysis as follows: under a reasonable
approximate mathematical assumption the method provides at least one point
in a sufficiently small neighborhood of the exact coefficient without 
\textit{a priori} knowledge of any point in this neighborhood. The accuracy of the
approximation or the distance between those points and the exact solution
depends on the error in the data and some parameters of the discretization.
We point out that the fact of the proximity of that point to the correct
coefficient, which was achieved without any \textit{a priori} knowledge of that small
neighborhood, is the main advantage of our globally convergent method over
locally convergent ones. Indeed, as soon as one knows a point in a
sufficiently small neighborhood of the true solution, one can refine it via
a locally convergent method, see, e.g.,~\cite[Chapter 4]{Beili2012}. The
latter, however, is outside of the scope of the current publication. We
refer to~\cite[Theorem 2.9.4]{Beili2012} for more
details about the definition of the global convergence as well as a rigorous
mathematical analysis of the global convergence of the method relying on an
approximate mathematical framework.

In previous works of the GCM summarized in~\cite{Beili2012}, the model
of hyperbolic wave type equations is considered and the time-domain problem
is converted into the pseudo-frequency domain problem via the Laplace
transform.  Since the Laplace transform
used in the method has an exponentially decaying kernel, one likely loses
some information in taking this transform of the (far field) measured data.
Therefore, we exploit the Fourier transform to improve the performance of the GCM and to extend its direct application
to multi-frequency data, which is common in applications to materials
characterization. This leads us to study a \textit{new} GCM in~\cite{Kliba2016}. More precisely, in the latter 
paper, we developed a GCM for solving the CIP for the Helmholtz equation
with multi-frequency data. The main difficulty in developing this new GCM is to work with complex-valued functions where the maximum principle, which plays an important role in the previous GCM~\cite{Beili2012}, is no longer applicable.

As a continuation of the work of~\cite{Kliba2016}, the goal of this paper is
to analyze the performance of our new GCM for multi-frequency
experimental raw data.  
There are some major new features of the present paper:
\begin{enumerate}
\item[i.]  The globally convergent approach
together with its advantages makes this work different from previously known
locally convergent methods.
\item[ii.]  This paper is the first one where we study the experimental
multi-frequency data for the new GCM of \cite{Kliba2016}. 
\item[iii.] For the multi-frequency raw data in this paper, we developed 
a new  data preprocessing procedure, which is discussed later in this section. This procedure is substantially different 
from that of the time domain data
in~\cite{Thanh2014}.
\item[iv.] Recall
that this is the CIP with a minimal amount of multi-frequency raw data (backscatter data associated
to a single incident wave) and we are not aware of any literature that addresses the numerical 
solution of this problem without using any advanced \textit{a priori} knowledge of the solution.
\end{enumerate}

 We also refer 
to~\cite{Kuzhu2012, Beili2012, Thanh2014} for our works on
time-domain measured data for the previous GCM.
 Furthermore, we remark that our measured data are
stable only on a small interval of frequencies surrounding the central
frequency of 3.1 GHz, see section~\ref{sect:centralFreq}. Therefore,
transforming these data into time domain by the inverse Fourier transform
may lead to a large inaccuracy in the time-dependent data obtained. That
means that time domain inversion techniques including the previous 
GCM~\cite{Beili2012} are not a good candidate for the inverse problem
under consideration.

Keeping in mind our desired application to the detection and identification
of mines and IEDs, we did not arrange any \textquotedblleft special"
conditions for our experiments, which would eliminate parasitic signals
scattered by some objects, which are outside of our interest. Thus, our data
have been collected in a regular room with the presence of office furniture,
computers, wifi signals, air conditioning, etc., see Figure~\ref{fi:1}.
Therefore, the measured data are contaminated by a significant amount of noise.
The latter is a major challenge for any inversion method. We refer to section~\ref%
{sect:dataProcess} for more details about the sources of the noise. We note
that conventional data denoising techniques are not helpful here due to the
rich structure of the measured data. We hence present in this paper a \emph{%
new} \emph{heuristic} data preprocessing procedure. This procedure aims to
make the raw data look somewhat similar to the corresponding computationally
simulated data. The latter is done basically via a sort of
\textquotedblleft distilling" signals reflected by the targets of our
interest from the rest of the measured signals. The preprocessed data are
then used as the input for the GCM. Our test results show that the GCM can
reconstruct with a good accuracy the dielectric constants of the scattering
objects as well as their geometric information such as location. In doing
so, we do not use any detailed \textit{a priori} knowledge on physical and
geometrical properties of the targets.

We also mention reconstruction algorithms of~\cite{Agalt2014, Novik2015, Kaban2004, Kaban2015} for
coefficient inverse problems, which obtain unknown coefficients without
using any advanced \textit{a priori} knowledge of their neighborhood. These
algorithms use data resulting from multiple measurement events, i.e., the
Dirichlet-to-Neumann map data. On the other hand, our GCM uses the data
corresponding to a single incident plane wave.

The paper is organized as follows. In section~\ref{sect:FP}, we state the
forward and inverse problems. Section~\ref{sect:FP4} is devoted
to a short description of the globally convergent method. We describe in
section~\ref{sect:imag} the data collection and the steps of data
preprocessing. Section~\ref{sect:result} contains a description of the
numerical implementation of the method and a summary of reconstruction
results. Summary of results can be found in section~\ref{sect:summary}.

\section{The forward and inverse problems}

\label{sect:FP}

\label{sec:2}
%
%

In this section we formulate the forward and inverse problems for scattering
of electromagnetic waves by a penetrable inhomogeneous medium in $\mathbb{R}%
^{3}$ in a certain range of the frequency $\omega $. Below ${\bf x}%
=\left( x,y,z\right) \in \mathbb{R}^{3}.$ We assume that the scattering
medium, which occupies a bounded domain $D\subset \mathbb{R}^{3}$, is
isotropic, non-magnetic and that it is characterized by the spatially
distributed dielectric constant $\varepsilon _{r}({\bf x})$, which is
also called the relative permittivity. 
Our forward problem consists in finding the function $%
u\left( {\bf x},k\right) $ from the following conditions:%
\begin{align}
\label{eq:Helm}
& \Delta u+k^{2}\varepsilon _{r}({\bf x})u=0,\quad {\bf x}\in \mathbb{R%
}^{3},  \\
& u=e^{ikz}+u_{\mathrm{sc}}, \\
 \label{eq:radiation}
& \lim_{r\rightarrow \infty }r\left( \frac{\partial u_{\mathrm{sc}}}{
\partial r}-iku_{\mathrm{sc}}\right) =0,\quad r=|{\bf x}|.
\end{align}
Here $k=\omega /c$ is the wavenumber with the speed of light $c$. Here $%
u\left( {\bf x},k\right) $ is one of components of the electric field.
The total wave field $u$ is the sum of the scattered field $u_{\mathrm{sc}}$
and the incident field $u_{\mathrm{inc}}=e^{ikz}$, propagating along the $z$%
-direction towards the scattering medium. We note that the scattered field $%
u_{\mathrm{sc}}$ satisfies the Sommerfeld radiation condition~%
\eqref{eq:radiation}, which guarantees that it is an outgoing wave.

\begin{remark}
 To justify the description of the propagation of 
electromagnetic waves by a the single Helmholtz equation \eqref{eq:Helm}
rather than by the full Maxwell's equation, we refer for instance to \cite[%
Chapter 13]{Born1999}. It is shown there that if the function $\varepsilon
_{r}({\bf x})$ varies slowly enough on the scales of the wavelength, then
the scattering problem for the Maxwell's equations can be approximated by
the scattering problem for the Helmholtz equation for a certain component of
the electric field. An additional justification comes out of the accuracy
of our reconstruction results for our experimental data, see Table \ref%
{tab:table2} in section 4.
\end{remark}

We assume that the function $\varepsilon _{r}({\bf x})$ is
frequency-independent, real-valued, $\varepsilon _{r}\in C^{1}(\mathbb{R}%
^{3})$ and that there exists a positive constant $\overline \varepsilon>1$ such that 
\begin{equation}
1\leq \varepsilon _{r}({\bf x})\leq \overline \varepsilon,\qquad \varepsilon _{r}({\bf x})=1\quad \text{for all }{\bf x}\in \mathbb{R}^{3}\setminus \overline{D}.
\label{eq:assum}
\end{equation}%
The latter condition in~\eqref{eq:assum} means that the medium outside of $D$
is homogeneous. It is well-known that the forward problem~\eqref{eq:Helm}--%
\eqref{eq:radiation} of finding the total field $u$, given $u_{\mathrm{inc}}$
and $\varepsilon _{r}$, is uniquely solvable for $u\in C^{2}(\mathbb{R}^{3})$%
, see~\cite[Chapter 8]{Colto1998}. Our assumption of the frequency
independence of $\varepsilon _{r}$ is justified by the fact that in our
inverse algorithm of \cite{Kliba2016} we actually work on a small interval
of wavenumbers $k$. Our direct experimental measurements of the dielectric
constants of targets have shown that they vary slowly with respect to $k$ on
small $k-$intervals.


We now formulate our inverse problem which is the main subject of the paper.
To this end let us define $\Omega \subset \mathbb{R}^{3}$ as a convex
bounded domain with its regular boundary $\partial \Omega$ such that $%
D\subset \Omega ,\partial D\cap \partial \Omega =\varnothing $. Let $%
\underline{k}$ and $\overline{k}$ be positive constants such that $%
\underline{k}<\overline{k}$. Let $\Gamma \subset \partial \Omega $ be the
part of the boundary $\partial \Omega $ which corresponds to the backscatter
side of $\Omega $. We consider the following inverse problem or
also the coefficient inverse problem:

\emph{Coefficient Inverse Problem.} Assume that we are given a
multi-frequency backscatter data $g({\bf x},k)$ defined by 
\begin{equation}
g({\bf x},k):=u({\bf x},k),\quad \text{for }{\bf x}\in \Gamma ,k\in
\lbrack \underline{k},\overline{k}],  \label{eq:CIP}
\end{equation}%
where $u({\bf x},k)$ is the total field in the forward problem~%
\eqref{eq:Helm}--\eqref{eq:radiation}. Determine the relative permittivity $%
\varepsilon _{r}({\bf x})$ for ${\bf x}\in \Omega $.

\begin{remark}
The mathematical justification of the globally convergent method~\cite%
{Kliba2016} has been done only for measurements on the entire boundary $%
\partial \Omega $. Thus, in our theoretical derivation we assume that the
function $g({\bf x},k)$ is given on $\partial \Omega$. In our numerical
implementation we would need to complement the backscatter data on the rest
of the boundary $\partial \Omega $. We show how this is done in section~\ref%
{sect:completion}.
\end{remark}

Since we consider only a single direction of the incident plane wave $u_{%
\mathrm{inc}}({\bf x},k)$, then this is an inverse problem with single
measurement data. Given such data, we emphasize that in our coefficient
inverse problem we aim to reconstruct the coefficient $\varepsilon _{r}(%
{\bf x})$ using a numerical algorithm. The question whether $\varepsilon
_{r}({\bf x})$ can be uniquely determined from the data is outside of the
scope of this paper. This is also one of the fundamental theoretical
questions in the field of inverse problems.

The first uniqueness result for multidimensional coefficient inverse
problems with single measurement data was established in~\cite{Bukhg1981},
where the authors introduced a method based on Carleman estimates. This
method has been extensively studied by a number of authors for uniqueness
theorems in inverse problems with a finite number of measurements. We refer
to~\cite{Yamam2009, Beili2012} and references therein for surveys of this
method and uniqueness results. However, the technique of \cite{Bukhg1981}
works only if a non-vanishing function $f({\bf x})$ stands in the right
hand side of~\eqref{eq:Helm}. Hence, we always assume below that uniqueness
result holds true.

We end up this section with the Lippmann-Schwinger equation formulation for
the forward problem~\eqref{eq:Helm}--\eqref{eq:radiation}. Let the wavenumber $k$ be fixed in~\eqref{eq:Helm}--\eqref{eq:radiation}. The Green's
function in free space for the forward problem is given by 
\begin{equation*}
\Phi _{k}({\bf x},{\bf y})=\frac{\exp {(ik|{\bf x}-{\bf y}|)}}{%
4\pi |{\bf x}-{\bf y}|},\quad {\bf x}\neq {\bf y}.
\end{equation*}%
It is well-known that, see~\cite[Chapter 8]{Colto1998}, the forward problem~%
\eqref{eq:Helm}--\eqref{eq:radiation} is equivalent to the
Lippmann-Schwinger equation 
\begin{equation}
u({\bf x})=e^{ikz}+k^{2}\int_{\Omega }\Phi _{k}({\bf x},{\bf y}%
)(\varepsilon _{r}({\bf y})-1)u({\bf y}) \d{{\bf y}},\quad {\bf x}%
\in \mathbb{R}^{3}.  \label{eq:LS}
\end{equation}%
We observe from~\eqref{eq:LS} that if we know the function $u({\bf x})$
in $\Omega $, then we can just extend it to ${\bf x}\in \mathbb{R}%
^{3}\setminus \overline{\Omega }$ by the integration in~\eqref{eq:LS}.
Hence, to solve~\eqref{eq:LS}, it is sufficient to find the function $u(%
{\bf x})$ only for points ${\bf x}\in \Omega $. This integral equation
plays an important role for the convergence analysis of \cite{Kliba2016} as
well as for the numerical algorithm of the globally convergent method, see
section~\ref{sect:Algorithm1}.

\section{The globally convergent numerical method}

\label{sect:FP4} We describe in this section the globally convergent method
together with its numerical algorithm. We essentially rely on the paper~\cite%
{Kliba2016} where a theoretical study for this version of the method has
been carried out for measurement data on the entire boundary $\partial
\Omega $. Therefore, even though we are given only experimental backscatter
data, we will assume that we have the complete data on $\partial \Omega$ for
our description of the method.

\subsection{Integro-differential equation formulation}

\label{sect:integral}

The first crucial step for the construction of the GCM is to formulate the
coefficient inverse problem as a nonlinear integro-differential equation.
This is different from locally convergent methods with iterative
optimization processes which typically rely on least-squares formulation. We
describe in this section how to derive this integro-differential equation
using the structure of the forward problem and a change of variables. To
follow the theory of~\cite{Kliba2016}, we assume below that numbers $%
\underline{k},\overline{k}$, which are given in~\eqref{eq:CIP}, are
sufficiently large and $k\in \lbrack \underline{k},\overline{k}]$. It was
shown in~\cite{Kliba2016} that $u({\bf x},k)\neq 0$ for ${\bf x}\in
\Omega $ as long as $k$ is sufficiently large.

It has been shown in~\cite{Kliba2016} that there exists a function $v(%
{\bf x},k)\in C^{2}(\overline{\Omega })$ such that $u({\bf x},k)=e^{v(%
{\bf x},k)}$. 
Substitution of $u({\bf x},k)=e^{v({\bf x},k)}$ in the Helmholtz
equation~\eqref{eq:Helm} and a simple calculation give 
\begin{equation}
\Delta v+\nabla v\cdot \nabla v=-k^{2}\varepsilon _{r}({\bf x}).
\label{eq:eqv}
\end{equation}%
From now on, for the sake of presentation, we write $a\cdot a$ as $a^{2}$
for a complex-valued vector $a\in \mathbb{C}^{3}$. We eliminate $\varepsilon
_{r}({\bf x})$ from~\eqref{eq:eqv} by the differentiation of both sides
of this equation with respect $k$, which is similar to the first step of the
method of \cite{Bukhg1981}. We obtain 
\begin{equation}
\Delta \partial _{k}v+2\nabla \partial _{k}v\cdot \nabla v=\frac{2\Delta
v+2(\nabla v)^{2}}{k}.  \label{eq:eqvk}
\end{equation}%
It is seen that if we can somehow find $v$ from the latter equation, then $%
\varepsilon _{r}$ can be computed via~\eqref{eq:eqv}. We observe in~%
\eqref{eq:eqvk} that $\partial _{k}v$ and $v$ are related as 
\begin{equation*}
v({\bf x},k)=-\int_{k}^{\overline{k}}\partial _{k}v({\bf x},s)\d s+v(%
{\bf x},\overline{k}).
\end{equation*}%
Defining 
\begin{equation}
q({\bf x},k):=\partial _{k}v({\bf x},k),\quad V({\bf x}):=v(\mathbf{%
x},\overline{k}),  \label{eq:defs}
\end{equation}%
and substituting in~\eqref{eq:eqvk}, we obtain a nonlinear
integro-differential equation for $q$ 
\begin{align}
& \frac{k}{2}\Delta q({\bf x},k)+k\nabla q({\bf x},k)\cdot \left(
-\int_{k}^{\overline{k}}\nabla q({\bf x},s) \d s+\nabla V({\bf x})\right) 
\notag  \label{eq:integralDif} \\
& =-\int_{k}^{\overline{k}}\Delta q({\bf x},s) \d s+\Delta V({\bf x}%
)+\left( -\int_{k}^{\overline{k}}\nabla q({\bf x},s) \d s+\nabla V({\bf x}%
)\right) ^{2}.
\end{align}%
Now we complete the integro-differential equation formulation for our
coefficient inverse problem by exploiting the data $g({\bf x},k)=u(%
{\bf x},k)$ given on the boundary $\partial \Omega $. Recall that $u(%
{\bf x},k)=e^{v({\bf x},k)},$ which deduces $\partial _{k}v({\bf x}%
,k)=\partial _{k}u({\bf x},k)/u({\bf x},k)$. We hence have a Dirichlet
boundary condition for $q({\bf x},k)$ on $\partial \Omega $ as 
\begin{equation}
q({\bf x},k)=\psi ({\bf x},k)\quad \text{on }\partial \Omega ,
\label{eq:boundarydata}
\end{equation}%
where $\psi ({\bf x},k):=\partial _{k}g({\bf x},k)/g({\bf x},k)$ on 
$\partial \Omega $.

We have obtained a nonlinear integro-differential equation~%
\eqref{eq:integralDif} for the function $q({\bf x},k)$ with the Dirichlet
boundary condition~\eqref{eq:boundarydata}. We call $V({\bf x})$ the tail
function. Both functions $q$ and $V$ in~\eqref{eq:integralDif} are unknown.
To solve our inverse problem, both these functions need to be approximated.
This is done by an iterative process, where we start the iterations from an
initial approximation $V_{0}({\bf x})$ for the tail function $V({\bf x}%
)$. We will describe later in section~\ref{sect:initialTail} how to find $%
V_{0}({\bf x})$, which is a crucial step in our method. Given $V_{0}$, we
find $q({\bf x},k)$ by solving~\eqref{eq:integralDif} and then compute $%
\varepsilon _{r}({\bf x})$ using~\eqref{eq:eqv}. This ends the first
iteration. The next iterations follow a similar procedure but they use an
updated tail function $V$ instead of its initial approximation $V_{0}$. Here
we observe that we need $\nabla V$ and $\Delta V$ instead of $V$ in our
computation. The gradient of the updated tail is computed by $\nabla V(%
{\bf x})=\nabla u({\bf x},\overline{k})/u({\bf x},\overline{k})$,
where $u({\bf x},\overline{k})$ is obtained by solving the
Lippman-Schwinger equation~\eqref{eq:LS} in $\Omega $ with $\varepsilon _{r}$
obtained from the previous iteration. Recall that $\Delta V=\mathrm{div}%
(\nabla V)$. This is an analog of the well known predictor-corrector
procedure, where updates for $V$ are predictors and updates for $q$ and $%
\varepsilon _{r}$ are correctors.

\subsection{Discretization with respect to the wavenumber}

\label{sec:5.4} To find $q$ and $V$ from~\eqref{eq:integralDif}--%
\eqref{eq:boundarydata} using an iterative process, we consider a
discretization of~\eqref{eq:integralDif}--\eqref{eq:boundarydata} with
respect to the wavenumber $k$. We divide the interval $[\underline{k},%
\overline{k}]$ into $N$ subintervals with the uniform step size $%
h=k_{n-1}-k_{n}$ as follows 
\begin{equation}
\underline{k}=k_{N}<k_{N-1}<...<k_{1}<k_{0}=\overline{k}.
\label{eq:divide_k}
\end{equation}%
We approximate the function $q({\bf x},k)$ as a piecewise constant
function with respect to $k\in \lbrack \underline{k},\overline{k}]$. 
More precisely, we assume that 
\begin{equation}
q({\bf x},k)=q_{n}({\bf x}),\quad \text{for }k\in \left[
k_{n},k_{n-1}\right) .  \label{eq:q_approx}
\end{equation}%
For $k\in \left[ k_{n},k_{n-1}\right) $, the latter approximation implies
that 
\begin{equation}
\int_{k}^{\overline{k}}q({\bf x},s)\d s=(k_{n-1}-k)q_{n}({\bf x}%
)+h\sum\limits_{j=0}^{n-1}q_{j}({\bf x}),  \label{eq:int_approx}
\end{equation}%
where $q_{0}({\bf x})=0$. Using~\eqref{eq:q_approx} and~%
\eqref{eq:int_approx} we rewrite problem~\eqref{eq:integralDif}--%
\eqref{eq:boundarydata} for $k\in \lbrack k_{n},k_{n-1})$ as 
\begin{align}
(k_{n-1}-k/2)\Delta q_{n}-F_{n}\cdot \nabla q_{n}-k_{n-1}(k_{n-1}-k)(\nabla
q_{n})^{2}& =G_{n}\quad \text{in }\Omega  \label{eq:discreteEq} \\
q_{n}& =\psi _{n}\quad \text{on }\partial \Omega ,  \notag
\end{align}%
where 
\begin{align*}
F_{n}& =(k+2(k_{n-1}-k))(h\sum\limits_{j=0}^{n-1}\nabla q_{j}+\nabla
V_{n-1}), \\
G_{n}& =(h\sum\limits_{j=0}^{n-1}\nabla
q_{j})^{2}-h\sum\limits_{j=0}^{n-1}\Delta q_{j}-2h\nabla V_{n-1}\cdot
\sum\limits_{j=0}^{n-1}\nabla q_{j}+\Delta V_{n-1}+(\nabla V_{n-1})^{2}.
\end{align*}%
We indicate here the dependence of the tail function $V:=V_{n}$ on $n$ since
we approximate $V$ iteratively as outlined in the end of section~\ref%
{sect:integral}. Now assuming that the step size $h$ is sufficiently small
and that $h\overline{k}\ll 1$, we ignore those terms in~\eqref{eq:discreteEq}%
, whose absolute values are $O(h)$ as $h\rightarrow 0$. We note that these
small terms include the nonlinear term $k_{n-1}(k_{n-1}-k)(\nabla q_{n})^{2}$
for $q_{n}$.

%

Even though the left hand side of equation~\eqref{eq:discreteEq} depends on $%
k$, it changes very little with respect to $k\in \left[ k_{n},k_{n-1}\right) 
$ since the interval $\left[ k_{n},k_{n-1}\right) $ is small. Hence, we
eliminate this $k-$dependence by integrating both sides of~%
\eqref{eq:discreteEq} with respect to $k\in \left( k_{n},k_{n-1}\right) $
and then dividing both sides of the resulting equation by $h$. We obtain the
Dirichlet boundary value problem 
\begin{align} 
\label{eq:bvp1}
\Delta q_{n}-\tilde{F}_{n}\cdot \nabla q_{n}& = \tilde{G}_{n}/k_{n-1} \quad 
\text{in }\Omega \\
 \label{eq:bvp2}
q_{n}& =\psi _{n}\quad \text{on }\partial \Omega,
\end{align}
where 
\begin{align*}
\tilde{F}_{n}& =(k_{n}/k_{n-1}+1)(h\sum\limits_{j=0}^{n-1}\nabla
q_{j}+\nabla V_{n-1}), \\
\tilde{G}_{n}& =-2h\sum\limits_{j=0}^{n-1}\Delta q_{j}-4h\nabla V_{n-1}\cdot
\sum\limits_{j=0}^{n-1}\nabla q_{j}+2\Delta V_{n-1}+2(\nabla V_{n-1})^{2}.
\end{align*}%
We refer to the paper~\cite{Kliba2016} for more details of the derivation of
the Dirichlet boundary value problem~\eqref{eq:bvp1}--\eqref{eq:bvp2} as
well as its solvability. Now we are ready to summarize the globally
convergent method as a computational algorithm in the next section.

\subsection{The algorithm}

\label{sect:Algorithm1}

In the algorithm below we use $n$ for the outer iterations and $i$ for the
inner iterations, where the latter play a role in updating the tails.


%
%
%
%
%
%
%
%

\begin{algorithm}[H]
\caption{Globally convergent algorithm}
\label{alg:globalconv}
\begin{algorithmic}[1]
\State{Given $\nabla V_{0}$, set $q_{0}:=0$}
\For{$n = 1, 2, \dots, N$}
\State{Set $q_{n, 0} := q_{n - 1}$ and $\nabla V_{n, 0} := \nabla V_{n - 1}$}
\For{$i=1, 2, \dots, I_N$}
\State{Find $q_{n,i}$ by solving the boundary value problem~\eqref{eq:bvp1}--\eqref{eq:bvp2}}
\State{Update 
$\nabla v_{n,i} :=- ( h\nabla q_{n,i}
+ h\sum_{j=0}^{n-1}\nabla q_{j}) +\nabla V_{n,i-1}$ in $\Omega$}
\State{Update $\eps_{rn,i}$ via~\eqref{eq:eqv}}
\State{Find $u_{n,i}(\x,\overline{k})$ by solving LS equation~\eqref{eq:LS} in $\Omega$ with $\eps_r :=\eps_{rn,i}$}
\State{Update  $\nabla V_{n,i}(\x) := \nabla u_{n,i}(\x,\overline{k})/ u_{n,i}(\x,\overline{k})$}
\EndFor
\State{Update $q_{n} :=q_{n,I_N}$, $\eps_{rn} :=\eps_{rn,I_N}$ and $\nabla V_n := \nabla V_{n,I_N}$}
\EndFor
\State{Choose $\eps_r$ by \emph{the-criterion-of-choice}}
\end{algorithmic}
\end{algorithm}



\begin{remark}
(i) The criterion for choosing the final result for $\varepsilon _{r}$ and
the stopping rules with respect to $n$ and $i$ in the iterations are
addressed in section~\ref{sect:stopping} of the numerical implementation. We
also refer to~\cite{Beili2012,Thanh2015} for stopping criteria developed for
the numerical verification for globally convergent methods. Recall that the
number of iterations is often considered as a regularization parameter in
the theory of ill-posed problems, see for instance~\cite{Beili2012}.

(ii) We refer to section~\ref{sect:details} for some details about the
truncation and smoothing that is used when updating $\varepsilon_{rn,i}$ via~%
\eqref{eq:eqv}.
\end{remark}

It can be seen from the algorithm that its main computational part of the
algorithm is to solve at each iteration the Dirichlet boundary value problem~%
\eqref{eq:bvp1}--\eqref{eq:bvp2} and the Lippmann-Schwinger equation~%
\eqref{eq:LS}. Recall that the well-posedness of these two problems have
been established, see~\cite{Kliba2016} for the boundary value problem and~%
\cite{Colto1998} for Lipmann-Schwinger equation. Some details about solving
these two problems can be found in the section~\ref{sect:details} of the
numerical implementation.

%
%

\subsection{The initial approximation $V_{0}({\bf x}) $ for the tail
function}

\label{sect:initialTail}

The initial approximation $V_{0}$ of the tail function plays an important
role in both numerical implementation and the convergence analysis. We refer
to~\cite{Kliba2016} for more details on how the convergence of the algorithm
depends on the initial tail function. Using an asymptotic expansion with
respect to $k$ for the solution $u({\bf x},k)$ to the forward problem~%
\eqref{eq:Helm}--\eqref{eq:radiation}, it has been shown in~\cite{Kliba2016}
that there exists a function $p({\bf x}) \in C^2(\overline{\Omega})$ such
that 
\begin{equation}
\log u({\bf x},k)=ik\,p({\bf x})(1+O(1/k)),\quad \text{as }%
k\rightarrow \infty .  \label{eq:asym}
\end{equation}%
Since $\overline{k}$ is sufficiently large, then dropping the term $O( 1/k) $
in~\eqref{eq:asym} for $k\geq $ $\overline{k},$ we obtain $v( {\bf x},k)
=\log u( {\bf x},k) =ikp ({\bf x}) $ for $k\geq $ $\overline{k}.$
Hence, we approximate the initial tail function as 
\begin{equation}  \label{eq:V_approx}
V( {\bf x}) =i\overline{k}p( {\bf x}).
\end{equation}
Since by~\eqref{eq:defs} $q( {\bf x},k) =\partial _{k}v( {\bf x} ,k)$,
then $q( {\bf x},\overline{k}) =ip( {\bf x})$. Substituting these
approximations in~\eqref{eq:integralDif} and~\eqref{eq:boundarydata} at $k=%
\overline{k}$, we obtain that the function $p( {\bf x}) $ is the solution
of the following Dirichlet boundary value problem for the Laplace equation 
\begin{align}  
\label{eq:p_Eq}
\Delta p({\bf x})& =0\quad \text{in }\Omega, \\
\label{eq:p_data}
p({\bf x})& =i\psi ({\bf x},\overline{k})\quad \text{on } \partial
\Omega.
\end{align}%
We refer to~\cite{Kliba2016} for more details on the derivation of the
latter boundary value problem as well as a discussion on its unique
solvability.

Therefore, for the initial approximation $V_{0}$ of the tail function $V$ we
take 
\begin{equation}  \label{eq:V0}
V_{0}({\bf x})=i\overline{k}\,p({\bf x})
\end{equation}%
where $p({\bf x})$ is computed by solving problem~\eqref{eq:p_Eq}--\eqref{eq:p_data}. 
It has been proved in~\cite{Kliba2016} that the accuracy
of this initial approximation depends only on the error in the boundary data 
\begin{equation*}
\Vert V_{0}({\bf x})-V({\bf x})\Vert _{C^{2,\alpha }(\overline{\Omega }%
)}\leq C\overline{k}\Vert \psi ({\bf x},\overline{k})-\psi _{\mathrm{exact%
}}({\bf x},\overline{k})\Vert _{C^{2,\alpha }(\partial \Omega )},
\end{equation*}%
where $C=C(\Omega )$ is a positive constant depending only on the domain $%
\Omega $ and $C^{2,\alpha}(\Omega)$ is the H\"{o}lder space. This means that
the error in the first tail function depends only on the error in the
boundary data. Therefore, we obtain a good approximation for the tail
function already on the zero iteration of our algorithm. However, our
numerical experience tells us that we need to do more iterations. The
approximation \eqref{eq:V_approx} and, as a result, the assumption~%
\eqref{eq:V0} form our reasonable approximate mathematical assumption
mentioned in Introduction. We note that~\eqref{eq:V_approx} and~\eqref{eq:V0}
are used only on the zero iteration of our algorithm and are not used for
other tail functions, which are obtained in the iterative process of our
algorithm.

With the choice~\eqref{eq:V0} of $V_{0}$, it was proved in~\cite{Kliba2016}
that the algorithm proposed in section~\ref{sect:Algorithm1} converges in
the sense that it provides a good approximation for the coefficient $%
\varepsilon _{r}({\bf x})$. The accuracy of this approximation depends
only on the noise in the backscatter data, discretization step size $h$ and
the domain $\Omega $. We note that no \textit{a priori} information about a small
neighborhood of the unknown coefficient is used here. Therefore we say that
our algorithm is globally convergent within the approximation framework~%
\eqref{eq:V_approx}.

\section{Experimental data and preprocessing}

\label{sect:imag} We describe in this section the experimental setup, how
the data were collected, and an important part of this paper: the data
preprocessing procedure.

\subsection{Data collection}

\begin{figure}[h!]
\centering
\subfloat[A photograph of our experimental
setup]{\includegraphics[width=5cm]{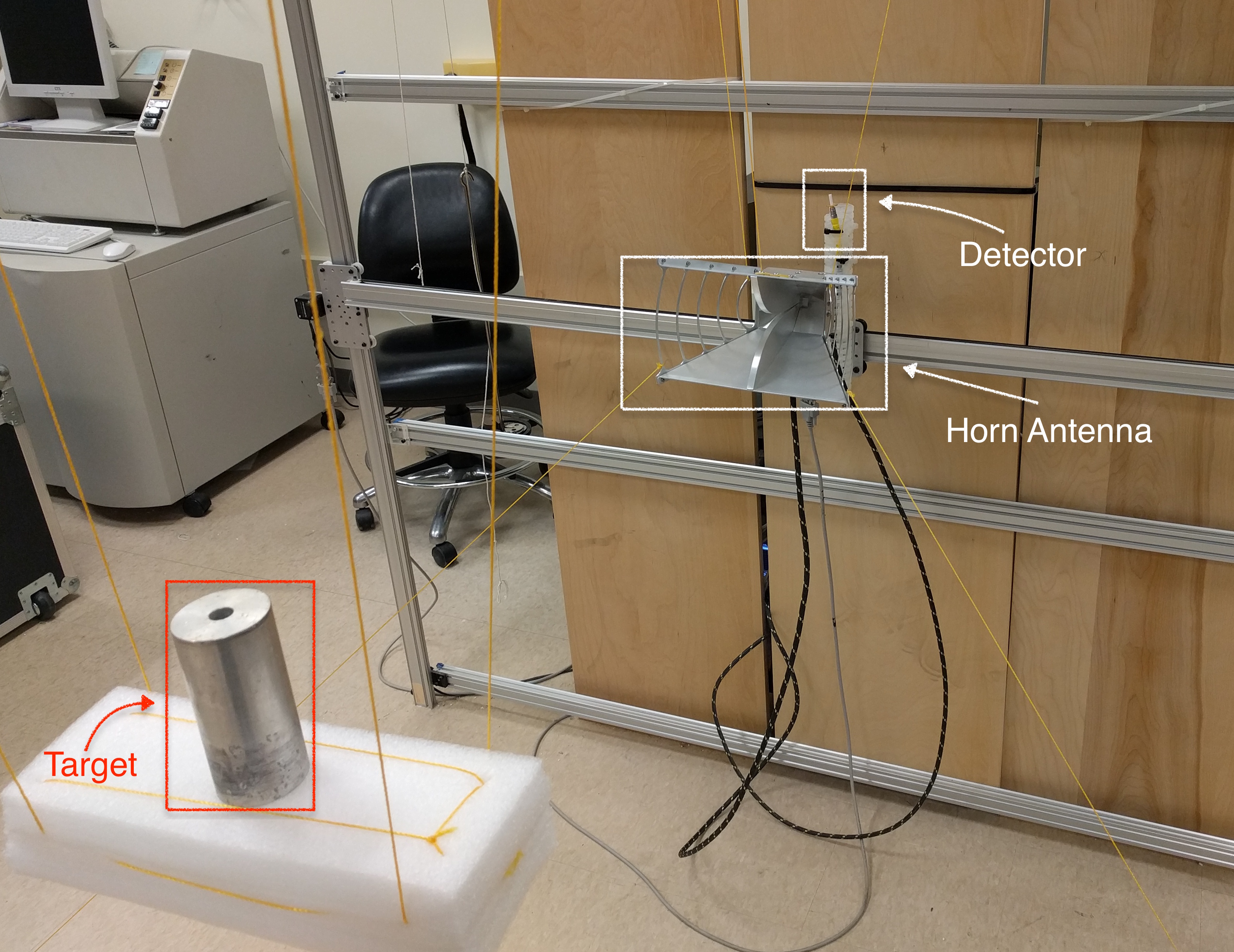}}\hspace{0.5cm} 
\subfloat[The
schematic diagram of our setup]{\includegraphics[width=6.5cm]{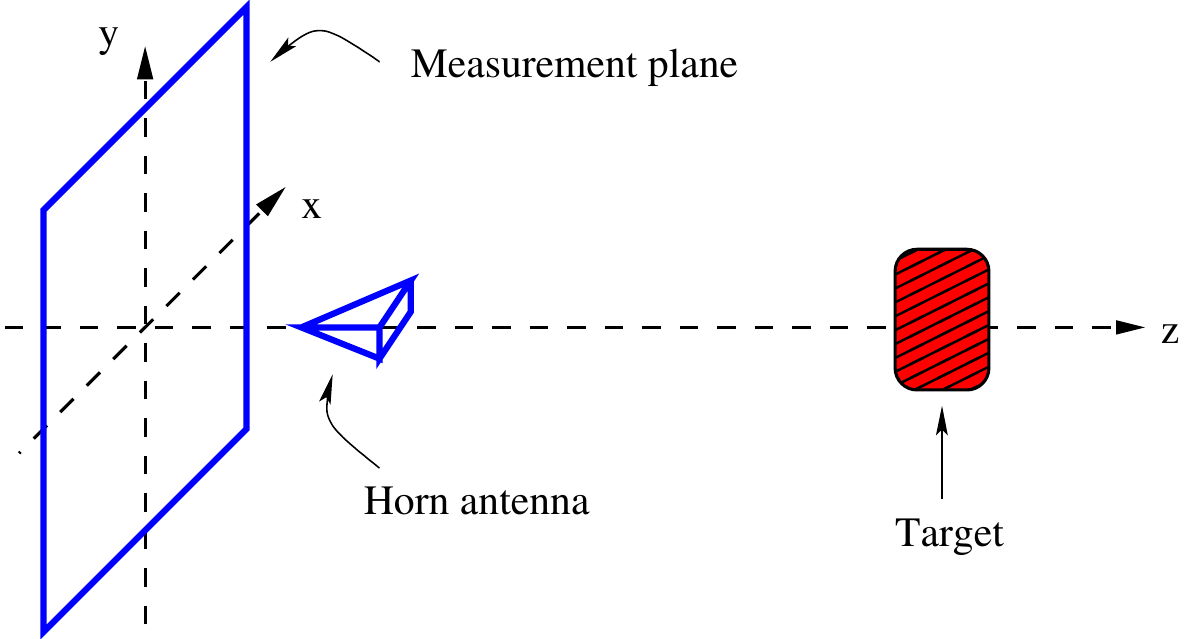}}
\caption{The experimental setup.}
\label{fi:1}
\end{figure}


The measurement data has been collected in a regular room, see Figure~\ref%
{fi:1}. As mentioned in Introduction, keeping in mind our desired
applications, we did not want to arrange a special anechoic chamber. On the
measurement plane, which is a 1 m by 1 m area, we define that the $x$-axis
and the $y$-axis are the horizontal and the vertical axis, respectively,
while the $z$-axis is perpendicular to the measurement plane. The direction
from the target, whose front face is positioned at $z=0$, to the measurement
plane is the negative direction of the $z$-axis.

The source and receiver in the experiments are both from the same device, a
2 port Rohde $\&$ Schwarz vector network analyzer (VNA), connected by
2-meter Megaphase RF Killer Bee test cables. The broadband horn antenna and
the collecting probe are both connected to the VNA. The antenna is
positioned at a distance of about 75 cm from the target and 18 cm from the
\textquotedblleft wooden wall" of Figure~\ref{fi:1}(a), where the probe
(i.e. detector) is located. For some technical reasons it was impossible to
place the antenna behind that \textquotedblleft wooden wall". The VNA sends
single frequency signals at frequencies ranging from 1 GHz to 10 GHz via the
horn antenna. Actually by design the signal is also emitted from the probe
at the same time. This probe technically collects both emitted and reflected
signals, more precisely their scattering parameters ($S$-parameters). 
However, in our setup, the reflection in the signal
to each device is calibrated when the detector is at the central position ($%
x=y=0$), that is, the internal reflections are removed at the position where
the signal is strongest. Nevertheless, the internal reflections are not calibrated at other locations of
the detector, where the signals are supposed to be weak. Therefore, 
the measured data are supposed to be as close to just the wave scattered 
by the unknown target as possible. Still, the latter approximation is only accurate 
at the central position, where the calibration of the internal reflections was done.
For the range of frequencies from 1 GHz to 10 GHz we collect the scattered
wave at 300 frequency points uniformly distributed over that range.

%

To obtain the backscatter data on the measurement plane for a single
incident wave, the experiment is repeated for different positions of the
probe on the measurement plane, which is the above mentioned wooden wall.
More precisely, the probe is uniformly moved over the scanning area with the
step size 2 cm, that means a total of 2550 data points. We scale to
dimensionless variables by considering ${\bf x}={\bf x}/(10\,\mathrm{cm%
})$ for the convenience in our numerical implementation. Therefore, from now
on when we see, for example $0.75$ of length, this means $7.5$ cm. Let $%
\mathbf{E}=(E_{x},E_{y},E_{z})$ be the electric field. The component $E_{y}$%
, which is the voltage, was incident upon the medium and the backscatter
signal of the same component was measured. Therefore, we denote $u({\bf x}%
,k):=E_{y}({\bf x},k)$, where the function $u({\bf x},k)$ is the above
discussed solution of the problem~\eqref{eq:Helm}--\eqref{eq:radiation}. As
it is stated in section~\ref{sect:centralFreq}, our central frequency is 3.1
GHz which corresponds to the wavelength $\lambda =9.67$ cm. This means that
the source/target distance was 7.75 $\lambda $. Since this distance is
sufficiently large in terms of $\lambda $, we have treated the incident wave
as a plane wave $E_{y,inc}=e^{ikz}$.

\subsection{Data preprocessing}

\label{sect:dataProcess}
\begin{figure}[h!]
\centering
\subfloat[Real part of measured
data]{\includegraphics[width=0.4\textwidth]{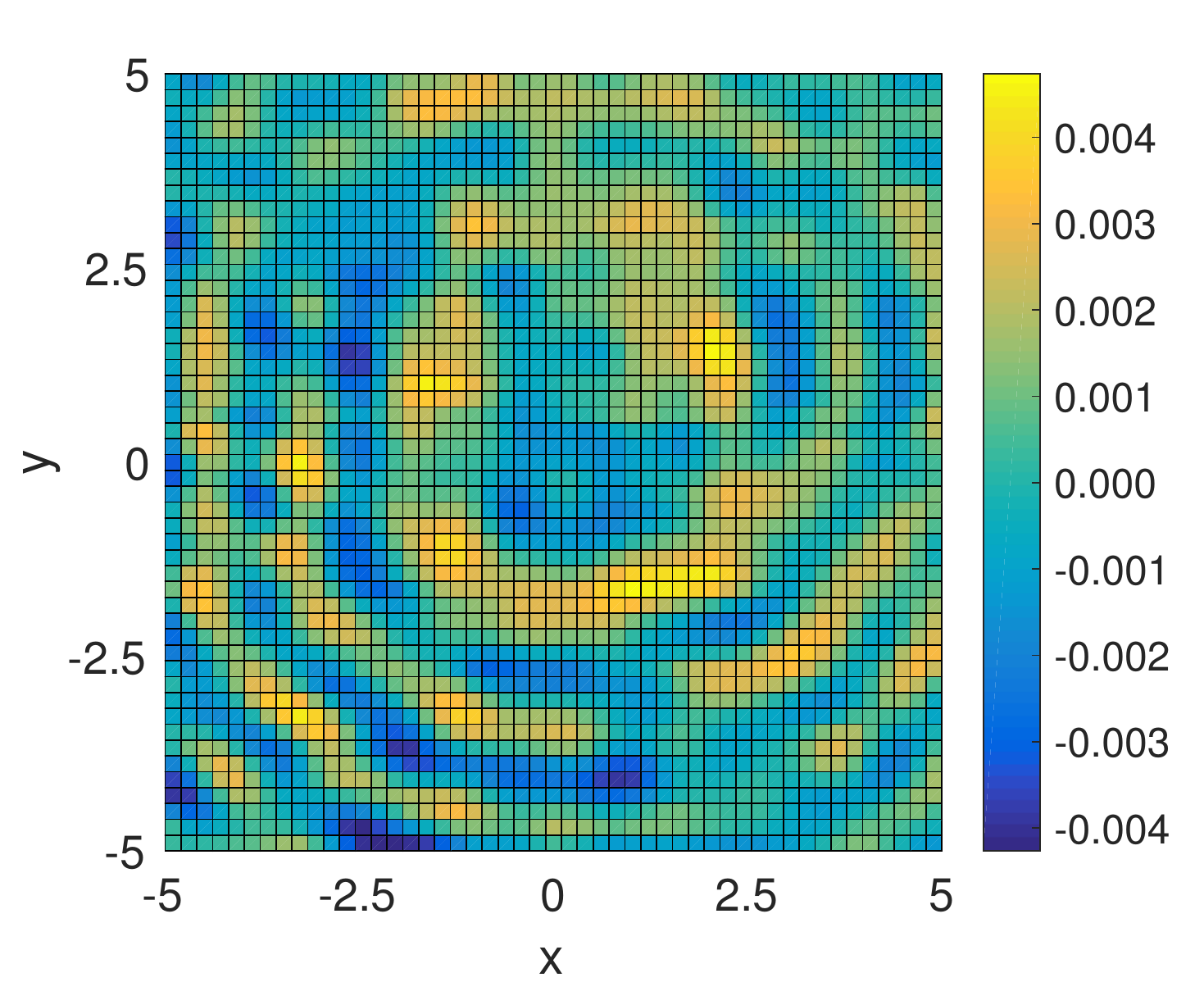}} \hspace{0.5cm} 
\subfloat[Real part of simulated
data]{\includegraphics[width=0.4\textwidth]{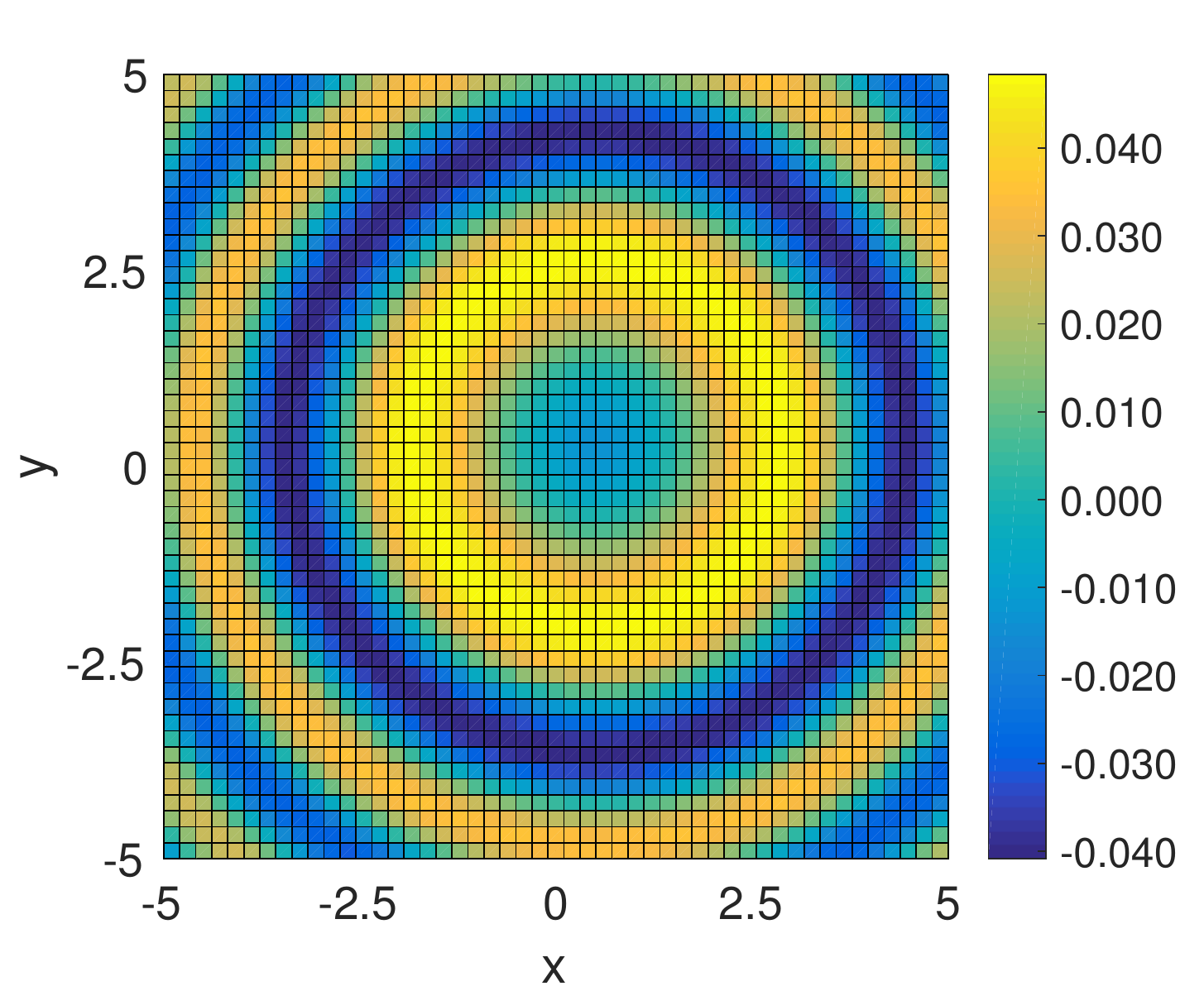}} 
\caption{The measured data (left) and the simulated data (right) on the
measurement plane ($xy$-plane) for $k=6.575$.}
\label{fi:2}
\end{figure}

Due to a significant amount of noise involved in the measured data,
there is a considerable mismatch between these data and the computationally
simulated data, see Figure~\ref{fi:2}. 
 Here are the main indicators for the
significant noise in the measured data. 
\begin{itemize}
\item[a.] The emitted signal propagates not only towards the target but
also  backwards to the measurement plane.

\item[b.] The measured signals are not only from the targets but also from other objects located
in the room. In addition, there were some metallic parts of our device placed behind the
measurement plane. So, the signals are scattered by
these parts and come back to the detectors. Furthermore, these signals are 
 affected by WiFi signals.

\item[c.]  The horn antenna is placed between the measurement plane and
the targets to be imaged, due to technical issues in the experimental setup. Hence, a part of the wave reflected by
these targets arrives back to the antenna. It is then scattered by the antenna and the latter signal propagates to the probe. This
causes some additional parasitic signals.

\item[d.] There is a significant difference in scales of magnitudes of experimental and computationally simulated
data. One of the reasons is that it is not clear how to  mathematically model the power of the emitted source.

\item[e.] Due to the instability of the
emitted signal, the backscatter signal is unstable. 
\end{itemize}

%
%

Therefore, to invert the measured data using the above numerical method, we
develop here a new heuristic data preprocessing procedure. The goal of this
procedure is twofold: (1) We sort of distill the signals reflected by
targets of our interest from signals reflected by other objects and (2) We
also make the experimental data look somehow similar to computationally
simulated data. It consists of three main steps:

\begin{enumerate}
\item[Step 1.] Data propagation. This step \textquotedblleft moves" the data
closer to the target. As a result, we obtain good estimates of $x,y$
coordinates of our targets. In addition, this step allows us to
reduce the computational domain $\Omega $.

\item[Step 2.] Selection of an interval of frequencies on which the
propagated data are stable. Our observation is that this is a quite narrow
interval.

\item[Step 3.] Truncation and calibration of propagated data.
\end{enumerate}

Next, the so preprocessed data are used as the input for the algorithm
described in section~\ref{sect:Algorithm1}.

\vspace{-0.0cm} 

\subsubsection{Data propagation}

This is one of the most important steps in data preprocessing. Given the
experimental data on the measurement plane, the data propagation process
aims to approximate these data on a plane which is closer to the unknown
target than the measurement plane. This process firstly reduces the
computational domain for our algorithm and secondly makes the scattered
field data to look more focused. In addition, we have observed in
our computations that the data propagation process helps us separate our
target signals from the unwanted signals, see the first paragraph of section \ref{sect:dataProcess}. This important observation leads to suitable truncations for
the data in the last step of the data preprocessing procedure.

%

The data propagation process is done using a well known method in optics
called the angular spectrum representation, see, e.g.,~\cite[Chapter 2]%
{Novot2012}. We note that this technique has been successfully exploited in~%
\cite{Thanh2015} for data preprocessing for the globally convergent method
of \cite{Beili2012}.

\begin{figure}[h]
\centering
\subfloat[Absolute value of  measured
data]{\includegraphics[width=0.4\textwidth]{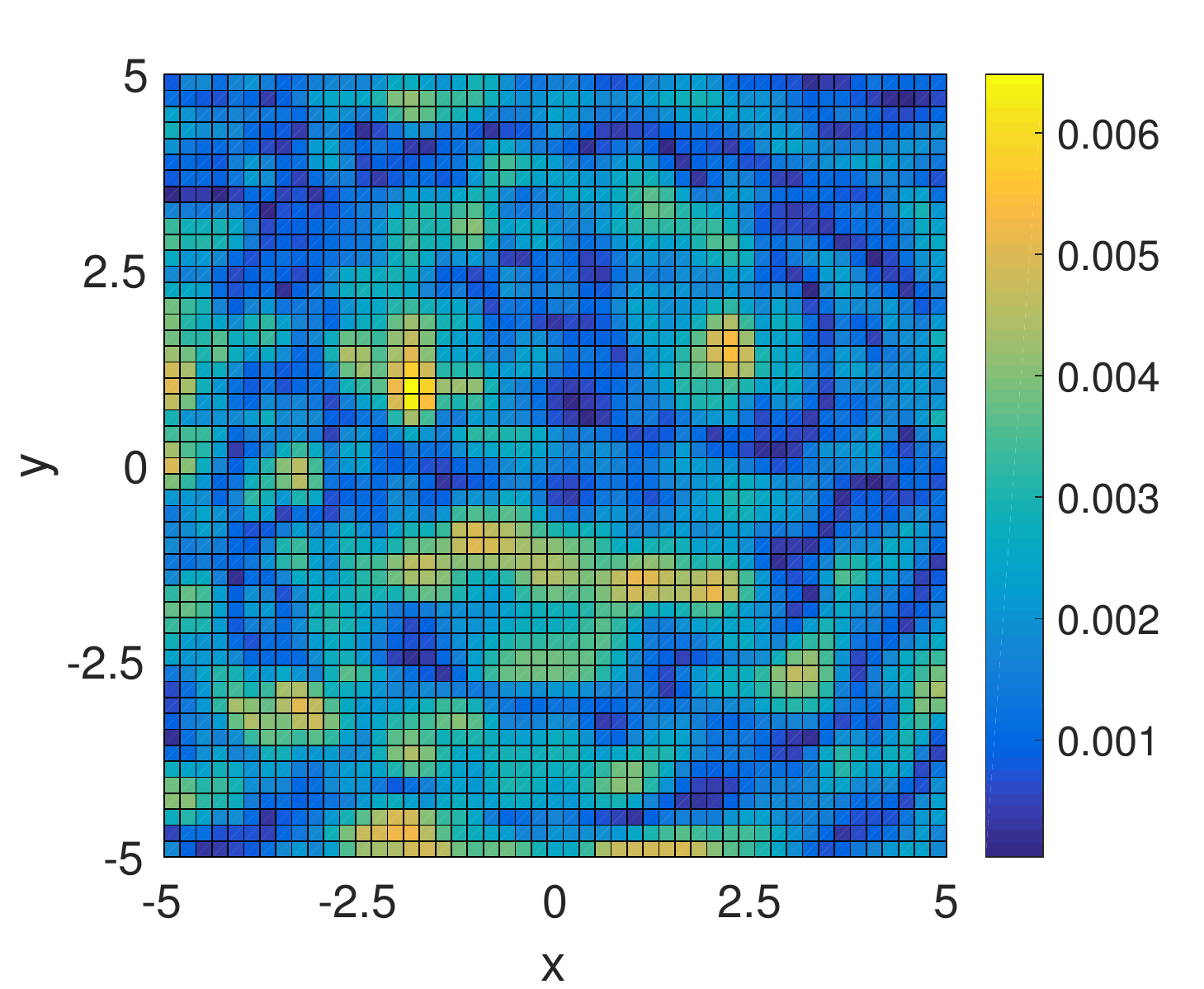}} \hspace{0.5cm} 
\subfloat[Absolute value of propagated
data]{\includegraphics[width=0.4\textwidth]{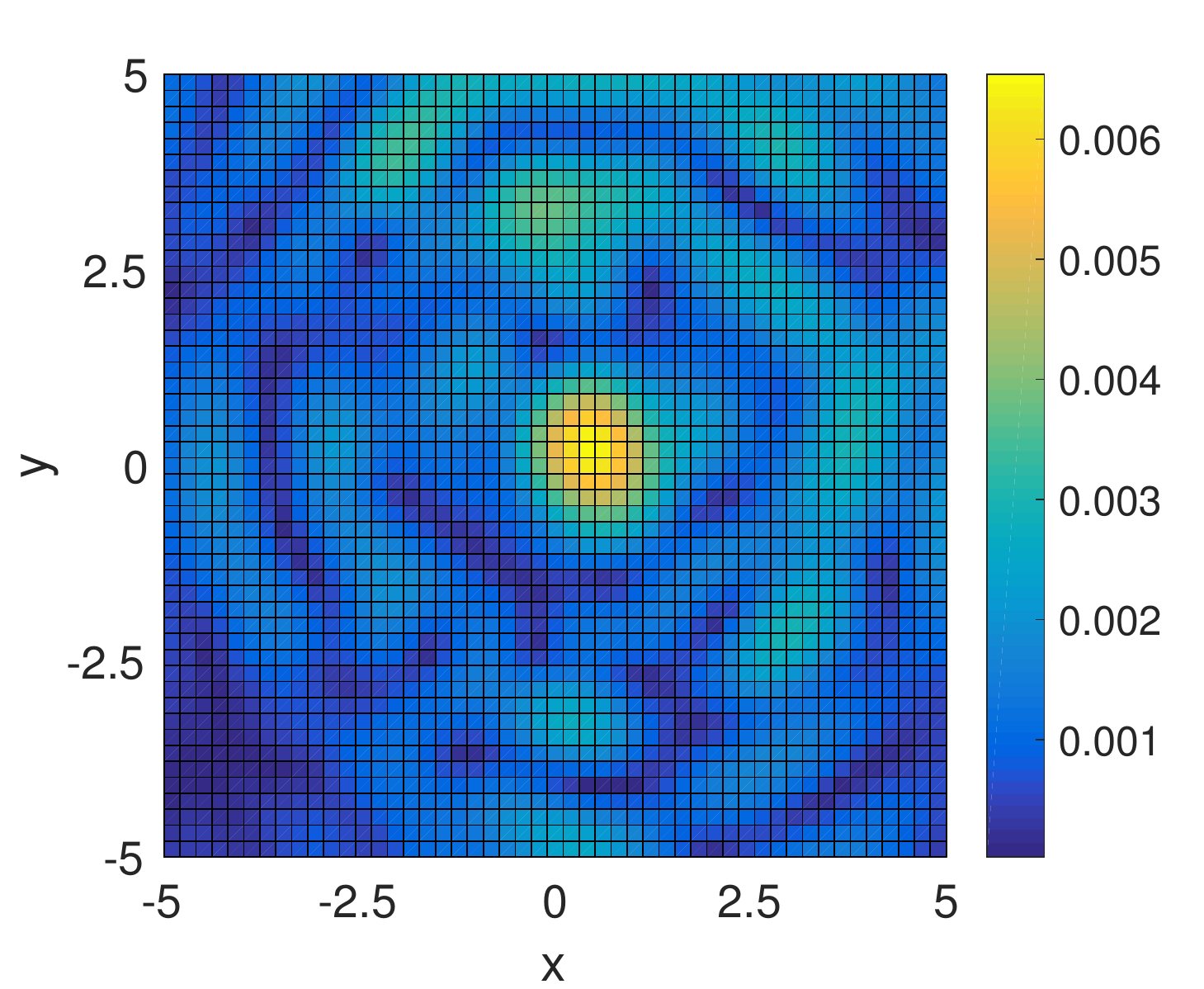}}
\caption{We respectively display in (a) and (b) the absolute value of the
measured and propagated data for $k=6.575$. }
\label{fi:3}
\end{figure}

We recall that the negative $z$-direction is from the target to the
measurement plane. Denote by $g({\bf x},k)$ the experimental data defined
on the measurement rectangle $P_{m}=\{{\bf x}:-5<x<5,-5<y<5,z=b\}$ of the
plane $\left\{ z=b\right\} $ and by $f({\bf x},k)$ its approximation on
the propagated rectangle $P_{p}=\{{\bf x}:-5<x<5,-5<y<5,z=a\}$ of the
plane $\left\{ z=a\right\} .$ We set both functions $g$ and $f$ to be equal
to zero on those parts of corresponding planes $\left\{ z=b\right\} $ and $%
\left\{ z=a\right\} $ which are outside of those rectangles. Abusing
notations a little bit, we call below sets $P_{m}$ and $P_{p}$
\textquotedblleft measurement plane" and \textquotedblleft propagated plane"
respectively. We have $b<a<0$ since $P_{p}$ is closer to the target than $%
P_{m}$. Using the method of the angular spectrum representation, we obtain 
\begin{equation}
f(x,y,a,k)=\frac{1}{2\pi }\iint_{k_{x}^{2}+k_{y}^{2}<k^{2}}\hat{g}%
(k_{x},k_{y},k)e^{-i[k_{x}x+k_{y}y-k_{z}(a-b)]} \d k_{x} \d k_{y},
\label{eq:propdata}
\end{equation}%
where 
\begin{equation*}
\hat{g}(k_{x},k_{y},k)=\frac{1}{2\pi }\iint_{\mathbb{R}%
^{2}}g(x,y,b,k)e^{i(k_{x}x+k_{y}y)} \d x \d y
\end{equation*}%
and $k_{z}=(k^{2}-k_{x}^{2}-k_{y}^{2})^{1/2}$. The idea of the method is
that we consider the Helmholtz equation satisfied by the total wave outside
the scattering medium ($\varepsilon _{r}({\bf x})=1$). Then the Fourier
transform with respect to $x$ and $y$ of the total field satisfies an
ordinary differential equation in the $z$-direction. Together with the
radiation condition and the boundary condition on the measurement plane $%
P_{m}$ one can solve the 1D problem and obtain the formula~%
\eqref{eq:propdata}. We refer to~\cite{Novot2012, Thanh2015} for more
details on the angular spectrum representation. In Figure~\ref{fi:3}(b) we
present the absolute value of the experimental data on the propagated plane 
\begin{equation}
P_{p}=\{{\bf x}:-5<x<5,-5<y<5,z=-0.75\},  \label{2}
\end{equation}%
computed by the formula~\eqref{eq:propdata} for the wavenumber $k=6.575$.
Comparing Figure~\ref{fi:3}(a) and Figure~\ref{fi:3}(b), we can easily see
that the focusing of the scattered field is considerably improved after data
propagation. From now on it is important to note that we will be only
interested in the propagated data instead of the measured one. \vspace{0cm} 

\subsubsection{Choosing a stable frequency interval}

\label{sect:centralFreq} 

\begin{figure}[!h]
\centering
{\ \includegraphics[width=0.4\textwidth]{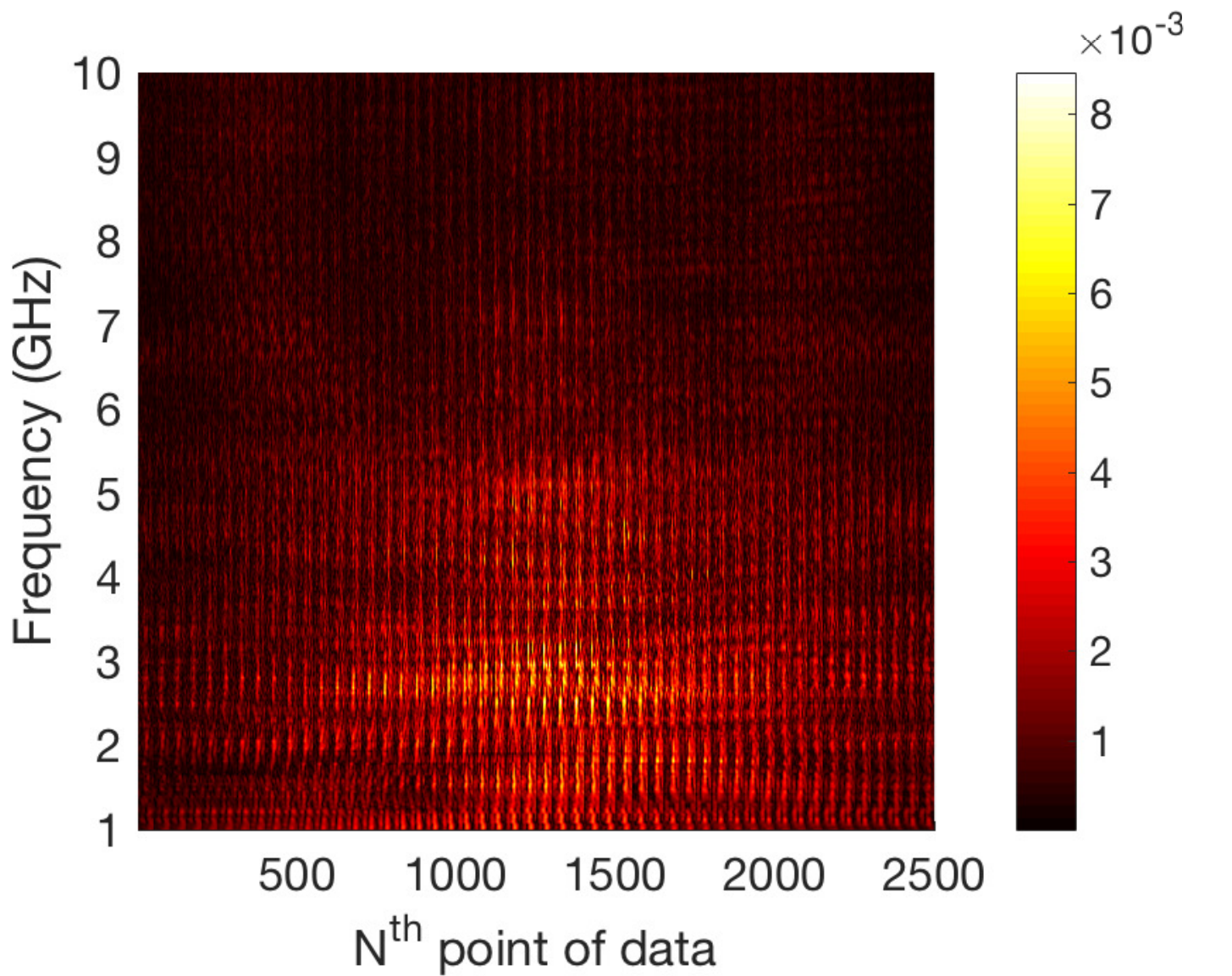}}
\caption{The absolute value of the experimental data after propagation for
all frequencies ranging from 1 GHz to 10 GHz. }
\label{fi:4}
\end{figure}

By manually plotting the absolute value of the propagated data for all 300
frequencies from 1 GHz to 10 GHz, we observed that the variation of the
focusing (as in Figure~\ref{fi:5}(b)) of the scattered field with respect to
the frequency is continuous only on a small interval of frequencies centered
at around 3.1 GHz. We call this interval a stable frequency interval, which
is essential for the input data for our numerical method. Equivalently, in
the numerical implementation this stable frequency interval determines the
corresponding small interval $[\underline{k},\overline{k}]$ for the wavenumbers that is used in our theoretical setting.
 We remark that this phenomenon does not exist for the computationally simulated data
since the propagated field focusing varies smoothly on any interval of
frequencies. The latter is one of the significant discrepancies between the
simulated data and the measured data. The stable frequency interval can also
be roughly estimated in Figure~\ref{fi:4}, where we plot the absolute value
of the propagated data, which is written as a matrix of $300\times 2550$.
The stable interval of frequencies from Figure~\ref{fi:4} is $[2.98,3.19]$
GHz, which implies that $[\underline{k},\overline{k}]=\left[ 6.25,6.70\right]
$. We observed that the lower and upper bounds of such intervals differ only
very slightly for the data from different targets considered in this paper.

From now on, by using $[\underline{k},\overline{k}]$ we understand that this
wavenumber interval corresponds to a small interval of frequencies around
3.1 GHz, which is determined from the procedure in this section.

\vspace{-0.0cm} 

\subsubsection{Truncation and calibration of propagated data}

\label{sect:truncation} 
\begin{figure}[!h]
\centering
\subfloat[Propagated data after
truncation]{\includegraphics[width=0.4\textwidth]{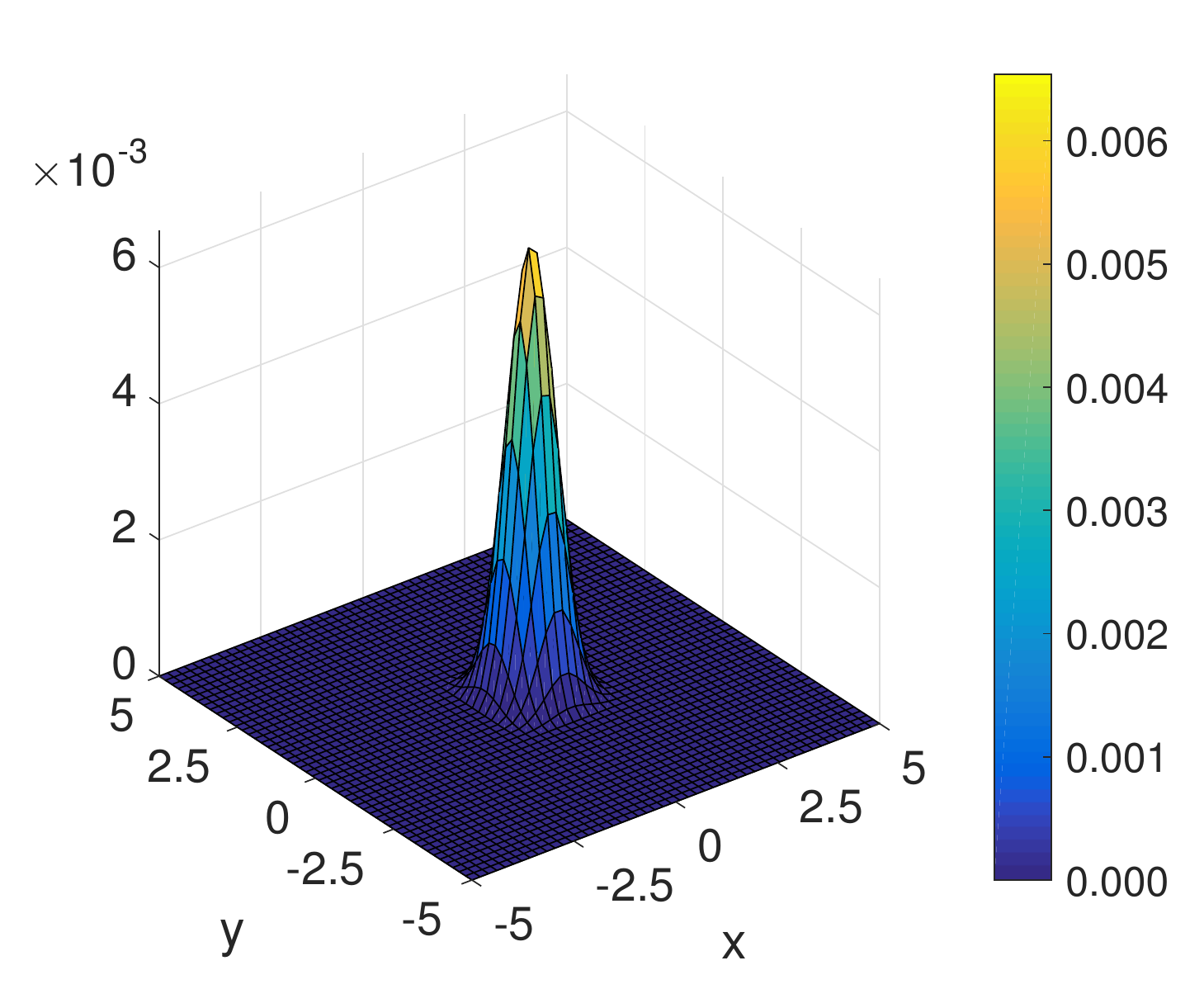}} \hspace{%
0.5cm} 
\subfloat[Propagated data in
simulation]{\includegraphics[width=0.4\textwidth]{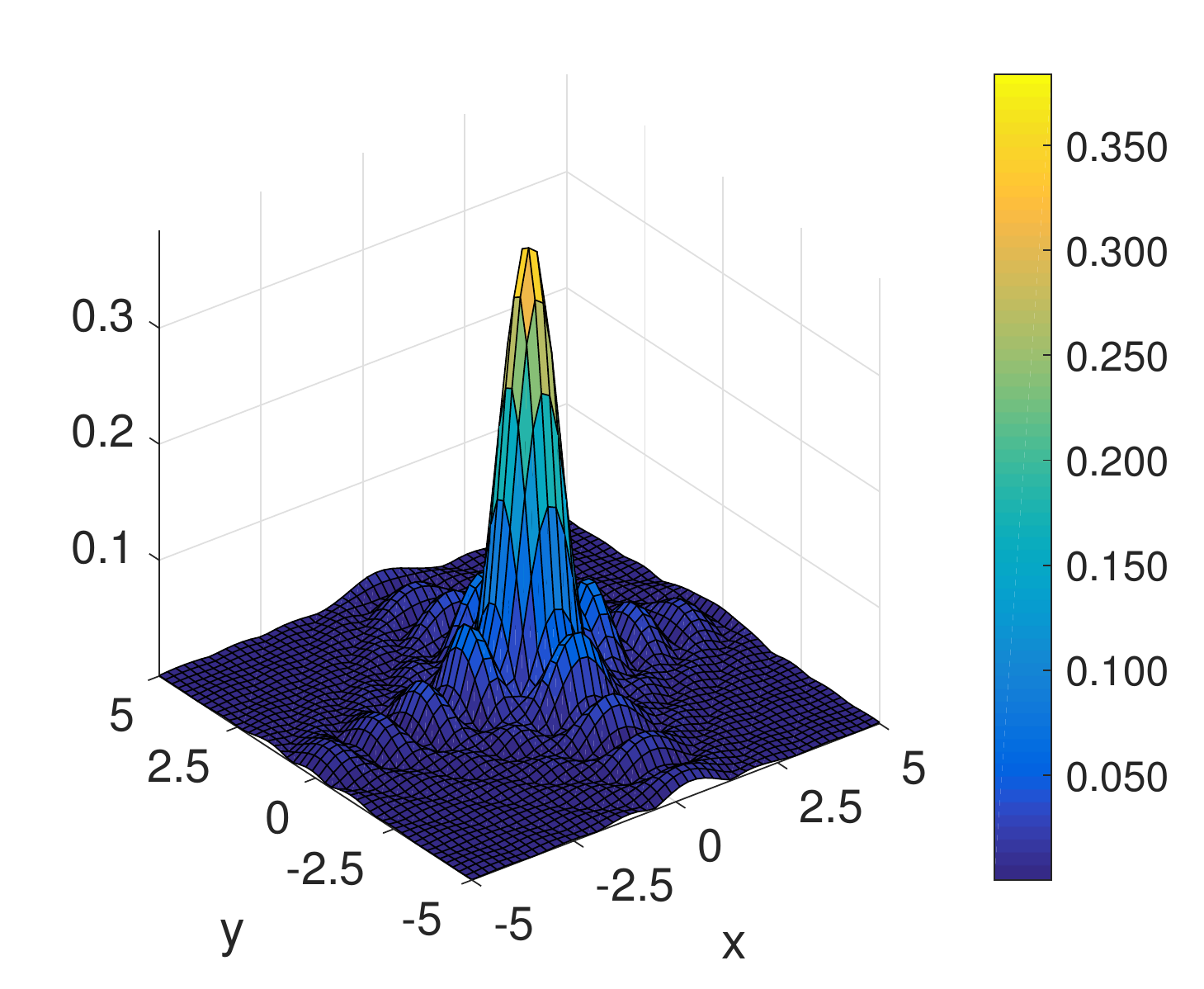}}
\caption{The absolute value of the propagated data after truncation and the
propagated data in simulation for $k=6.575$.}
\label{fi:5}
\end{figure}

We know from section~\ref{sect:centralFreq} that the
focusing of the propagated data changes continuously on $[\underline{k},%
\overline{k}]$. However, the behavior of the absolute value of the
propagated data on this interval of wavenumbers is still not close enough
to those in simulations. Indeed, we observed that after data propagation the
magnitude of the simulated data has a clear peak at the $xy$-location of the
target for all frequencies, see Figure~\ref{fi:5}(b). For the propagated
data, although its magnitude always attains the global maximum at the $xy$%
-location of the target for all $k\in \lbrack \underline{k},\overline{k}]$,
it does not really have a clear peak as in simulations due to local maxima
at other regions of the propagated plane, see Figure~\ref{fi:3}(b). However,
those latter local maxima are not located at the same locations for
different frequencies. On the other hand the global maximum is located at
the same place for all $k\in \lbrack \underline{k},\overline{k}]$. From this
important observation we guess that the $xy$-location of the target should
be actually in some small neighborhood of the point where the propagated
data attains its global maximum.


To make the propagated data look more similar to those in simulations we
truncate those by taking the points for which its magnitude is greater than
or equal to $80\%$ of its maximum value, and setting the rest to be zero. In
other words, for each $k\in \lbrack \underline{k},\overline{k}]$ we replace
the function $f(x,y,a,k)$ with the function $\widetilde{f}(x,y,a,k)$ where 
\begin{equation}
\widetilde{f}(x,y,a,k)=%
\begin{cases}
f(x,y,a,k), & \text{ if }|f(x,y,a,k)|\geq 0.8\max_{x,y}|f(x,y,a,k)|, \\ 
0, & \text{ otherwise.}%
\end{cases}
\label{eq:truncation}
\end{equation}%
After that we smooth the function $\widetilde{f}(x,y,a,k)$ by a Gaussian
filter and obtain the function $\widetilde{f}_{smth}(x,y,a,k)$. More
precisely, we use the smoothing function \texttt{smooth3} in Matlab with the
Gaussian option. This is done for all values of $k\in \lbrack \underline{k},%
\overline{k}]$. We display in Figure~\ref{fi:5} the propagated data after
this procedure for $k=6.575$. It can be seen that the behavior of the
absolute value of these processed data looks more similar to that of
simulated data. The results are similar for other wavenumbers $k$ in $[%
\underline{k},\overline{k}]$.

We are now at the last step of data preprocessing called data calibration.
This process is to scale the amplitude of the measured data to the same
scaling as in our simulations, since the experimental data and simulated
data usually have quite different amplitudes. For example, in Figure~\ref%
{fi:5}, we can see that the maximal value of the absolute value of measured
data is about 0.0065, while this value is about 0.4 for the simulated data.
To do the scaling, we need the so-called calibration factor. For estimating
the calibration factor, we use the measured data of a single calibrating
object whose location, shape, size, and dielectric constant are known. The
word ``single" is important here to ensure that the calibration procedure is
unbiased, which means that we are supposed to know such information for only
the calibrating object.

For $k\in [\underline{k},\overline{k}]$, where $\underline{k}$ and $%
\overline{k}$ are chosen from the process of choosing a stable frequency
interval of section~\ref{sect:centralFreq}, let $U_{cal.obj}^{sim}({\bf x}%
,k)$ and $U_{cal.obj}^{exp}({\bf x},k)$ be the propagated data of the
simulated and measured scattered field for the calibrating object,
respectively. We define our calibration factor as 
\begin{equation*}
\quad A(k)=\frac{\max {|U_{cal.obj}^{sim}({\bf x},k)|}}{\max {\
|U_{cal.obj}^{exp}({\bf x},k)|}}.
\end{equation*}
It can be seen that the amplitude of $A(k)U_{cal.obj}^{exp}({\bf x},k)$
has the same scale to that of the simulated data. We use the calibration
factor $A(k)$ for measured data of all other objects. More precisely, given
the propagated and truncated scattered field data $\widetilde{f}_{smth}(
x,y,a,k)$ obtained from~\eqref{eq:truncation}, we consider $A(k) \widetilde{f%
}_{smth}( x,y,a,k) + e^{ikz}$ as the input for our algorithm in section~\ref%
{sect:Algorithm1}. We emphasize that this calibration factor is computed
using the geometrical and physical information of only a single object: the
calibrating object. \vspace{-0.1cm} 

\section{Numerical implementation and reconstruction results}

\label{sect:result} We describe in this section some details of the
numerical implementation and present the reconstruction results for our
experimental data using the globally convergent algorithm. We have collected
experimental data for five objects. So, we test here these five data sets.
Each data set corresponds to a single scattering object, numbered from 1 to
5, see Table~\ref{tab:table1}. Object 1 (a piece of yellow pine) is chosen
as our calibrating object. \vspace{-0.0cm} 
\begin{table}[tbph]
\caption{Scattering objects}
\label{tab:table1}\centering
\begin{tabular}{|c|c|}
\hline
Object ID & Name of the object \\ \hline
1 & A piece of yellow pine \\ \hline
2 & A piece of wet wood \\ \hline
3 & A geode \\ \hline
4 & A tennis ball \\ \hline
5 & A baseball \\ \hline
\end{tabular}%
\end{table}

\subsection{Some computational details}

\label{sect:details} In our numerical implementation the front face of the
target is positioned at $\{z=0\}$. We consider the domain $\Omega$ as 
\begin{align*}
\Omega =(-2.5,2.5)\times (-2.5,2.5)\times (-0.75,4.25).
\end{align*}
Here we propagate all the measured data to the rectangle $P_{p}$ located on
the plane $\{z=-0.75\}$ and search for the unknown targets in the range $%
(-0.75,1)$ of the $z$-direction. This search range is motivated by both:
similar ranges that have been used in~\cite{Thanh2014,Thanh2015} and our
desired applications for detection of hidden explosives. Indeed, $z\in
(-0.75,1)$ means that the linear size of the target in the $z$-direction
does not exceed 17.5 cm. It is well-known that typical linear size of
antipersonnel mines and IEDs are between 5 cm and 15 cm. 
We note that the larger range $(-0.75,4.25)$ is considered for the goal of
completing the backscatter data in the next section. Also in the $xy$-plane
we restrict ourselves to the smaller area $(-2.5,2.5)^{2}$ than the original
measurement one $(-5,5)^{2}$. The main reason comes from the observation
that the truncated data is zero outside of the area $(-2.5,2.5)^{2}$, see
Figure~\ref{fi:5}, and does not contribute to the reconstruction process.
This restriction helps us reduce our computational domain.

To solve the Lippmann-Schwinger integral equation~\eqref{eq:LS}, we use a
spectral method developed in~\cite{Lechl2014, Vaini2000}. This method
exploits a periodization technique of the integral equation introduced by
G. Vainikko, which enables the use of the fast Fourier transform and a
simple numerical implementation. We use this method to find the
solution in a box containing the support of the function $\varepsilon_r -1$.
The solution in the complement of the box in $\Omega$ is then computed by
the integration in~\eqref{eq:LS}. We solve the boundary value problem~%
\eqref{eq:bvp1}--\eqref{eq:bvp2} by a finite element method. More precisely,
we use FreeFem++~\cite{Hecht2012}, a standard software designed with a focus
on solving partial differential equations using finite element methods. We
refer to \emph{www.freefem.org} for more information about FreeFem++.

We observed from data preprocessing that the location of a target in the $xy$%
-plane can be roughly estimated from the propagated data. Indeed, we define $%
\Omega _{T}$ as 
\begin{equation*}
\Omega _{T}=\left\{ (x,y):|\widetilde{f}_{smth}(x,y,-0.75,\widetilde{k}%
)|>0.7\max |\widetilde{f}_{smth}(x,y,-0.75,\widetilde{k})|\right\} \subset
P_{p},
\end{equation*}%
where $\widetilde{f}_{smth}(x,y,-0.75,\widetilde{k})$ is the propagated data
after the truncation procedure in section~\ref{sect:truncation} and $%
\widetilde{k}$ is the wavenumber corresponding to the frequency 3.1 GHz.
Note that for each target $|\widetilde{f}_{smth}(x,y,-0.75,k)|$ has a
positive peak whose location is the same for all $k\in \lbrack \underline{k},%
\overline{k}]$, see Figure~\ref{fi:5}. The truncation value 0.7 was chosen
based on trial-and-error tests on simulated and calibrating targets. We
observed that $\Omega _{T}$ provides a reasonable approximation for the $xy$%
-location of a target. The same truncation was applied to all targets.
Hence, it is not biased. Using $\Omega _{T}$ and the assumption that we seek
for the target in $(-0.75,1)$ of the $z$-direction, the coefficient $%
\varepsilon _{rn,i}$ is truncated as follows 
\begin{equation*}
\varepsilon _{rn,i}({\bf x}):=%
\begin{cases}
\max \left( |\varepsilon _{rn,i}({\bf x})|,1\right) , & {\bf x}\in
\Omega _{T}\times (-0.75,1), \\ 
1, & \text{elsewhere}.%
\end{cases}%
\end{equation*}%
After that it is smoothed by \texttt{smooth3} in Matlab with the Gaussian
option as in data truncation. The so obtained $\varepsilon _{rn,i}$ is then
used for solving the Lippmann-Schwinger equation in Algorithm~\ref%
{alg:globalconv}.

Another important detail we would like to mention is a modification for
approximating the data $\psi ({\bf x},k)=\partial _{k}g({\bf x},k)/g(%
{\bf x},k)$ on $\Gamma =\{{\bf x}\in \partial \Omega :z=-0.75\}$,
where $g({\bf x},k)=A(k)\widetilde{f}_{smth}(x,y,-0.75,k)+e^{-0.75ik}$.
First we calculated $\partial _{k}g({\bf x},k)$ using finite difference.
The function $g({\bf x},k)$ is supposed to be as close to the total field
as possible and $\psi ({\bf x},k)$ should be smooth with respect to $k$.
However, the variation of $\psi ({\bf x},k)$ with respect to $k$ in the
small interval~$[\underline{k},\overline{k}]$ is actually not small, since
the behavior of $g({\bf x},k)$ with respect to $k$ is not similar to that
of the corresponding total field in simulations. The latter is one of
discrepancies between simulated and propagated experimental data. We observe
from simulations that, for each $k\in [\underline{k},\overline{k}]$, the
absolute value of the function $\psi ({\bf x},k)=\partial _{k}g({\bf x}%
,k)/g({\bf x},k)$ on the propagated plane typically has a global minimum
at the $xy$-location of the scatterer. Motivated by this observation, we
approximate $\psi ({\bf x},k)$ by $\widetilde{\psi }({\bf x}%
,k)=\partial _{k}g({\bf x},k)/g({\bf x} ,k^{\ast })$, where for each $%
k $ we choose such a value $k^{\ast } \in [\underline{k},\overline{k}]$ for
which $|\widetilde{\psi }({\bf x},k)|$ has a global minimum. However, we
also observed that the choice of $k^*$ is independent of $k$. By doing so,
we see that the variation of $\widetilde{\psi }({\bf x},k)$ in $k$
becomes small enough, when compared with those in simulations to be used in
the algorithm.

\subsection{Completing the backscatter data}

\label{sect:completion}

We recall that only the backscatter signals were measured in our
experiments. This means that after the above data preprocessing procedure,
the function $g({\bf x},k)$ was known only on the side $\Gamma =\{\mathbf{%
x}\in \partial \Omega :z=-0.75\}$ of the domain $\Omega $. As in~\cite%
{Thanh2014, Kliba2016}, we replace the missing data on the other parts of $%
\Omega $ by the corresponding solution of the forward problem in the
homogeneous medium, where $\varepsilon_r({\bf x}) =1$. In other words, in
our computation, the function $g( {\bf x},k)$ is extended on the entire
boundary $\partial \Omega $ as 
\begin{equation}  \label{eq:completion}
g({\bf x},k):= 
\begin{cases}
g({\bf x},k), & {\bf x}\in \Gamma , \\ 
e^{ikz}, & {\bf x}\in \partial \Omega \setminus \Gamma .%
\end{cases}%
\end{equation}
We remark that data completion methods are widely used for inverse problems
with incomplete data. The data completion~\eqref{eq:completion} is a
heuristic technique relying on the successful experiences of our group when
working with globally convergent methods for experimental data, see~\cite%
{Thanh2014, Beili2012a}. Other data completion methods may also be applied.

\subsection{The stopping criteria and choosing the final result}

\label{sect:stopping}

Now we introduce the rules for stopping the iterations and choosing the
final result for our Algorithm~\ref{alg:globalconv} relying on the content
of the convergence theorem of~\cite{Kliba2016} and trial-and-error testing
for simulated data. Indeed, the global convergence theorem 7.1 of~\cite%
{Kliba2016} only claims that functions $\varepsilon_{rn,i}$ are located in a
sufficiently small neighborhood of the true solution if the number of
iterations is not too large. But this theorem does not claim that these
functions tend to the exact solution, also see theorem 2.9.4 in~\cite%
{Beili2012} and theorem 5.1 in~\cite{Beili2012a}. This implies that the
stopping criteria should be chosen computationally. 

We have two stopping criteria: one for the inner iterations and the second
one for the outer iterations. Before stating those criteria precisely we
need some definitions. Denote by $e_{n,i}$ the relative error between the
two computed coefficients corresponding to two consecutive inner iterations
of the $n$-th outer iteration. That means 
\begin{equation*}
e_{n,i} = \frac{\Vert \varepsilon _{rn,i}-\varepsilon _{rn,i-1}\Vert
_{L^{2}(\Omega )}}{\Vert \varepsilon _{rn,i-1}\Vert _{L^{2}(\Omega )}},
\quad \text{for } i = 2, 3,\dots
\end{equation*}
Consider the $n$-th and $(n+1)$-th outer iterations which contain $I_1$ and $%
I_2$ inner iterations, respectively. We define the sequence of relative
errors associated to these two outer iterations as 
\begin{equation}  \label{eq:sequence}
e_{n,2}, \dots , e_{n,I_1}, \widetilde{e}_{n+1,1}, e_{n+1,2}, \dots ,
e_{n+1,I_2},
\end{equation}
where 
\begin{equation*}
\widetilde{e}_{n+1,1} = \frac{\Vert \varepsilon _{rn+1,1}-\varepsilon
_{rn,I_1}\Vert _{L^{2}(\Omega )}}{\Vert \varepsilon _{rn,I_1}\Vert
_{L^{2}(\Omega )}}.
\end{equation*}

The inner iterations with respect to $i$ of the $n$-th outer iteration in
the algorithm~\ref{alg:globalconv} are stopped when either $e_{n,2} <
10^{-6} $ or $i=3$. Note that this rule is similar to the one used for
``Test 2" in~\cite{Thanh2014}, where the maximal number of inner iterations
is set up to be 5. During our numerical experiments we observed that the
reconstruction results are essentially the same when we use either 3 or 5
for the maximal number of inner iterations.

Concerning the outer iterations with respect to $n$ in Algorithm~\ref%
{alg:globalconv}, it can be seen from the stopping rule for the inner
iterations that each outer iteration consists of at least 2 and at most 3
inner iterations. Equivalently, the error sequence~\eqref{eq:sequence} has
at least 3 and at most 5 elements. We stop the outer iterations if there are
two consecutive outer iterations for which their error sequence~%
\eqref{eq:sequence} has three consecutive elements less than or equal to $%
5\times 10^{-4}$. 
We emphasize again that we have no rigorous justification for these stopping
rules; they rely on the content of the convergence theorem~\cite{Kliba2016}
and trial-and-error testing for simulated data.

We choose the final result for $\varepsilon _{r}({\bf x})$ by taking the
average of its approximations $\varepsilon _{rn,i}({\bf x})$
corresponding to the relative errors in~\eqref{eq:sequence} that meet the
stopping criterion for outer iterations. The computed dielectric constant is
determined as the maximal value of the computed $\varepsilon _{r}({\bf x}%
) $. We have observed in our numerical studies that we need no more than
five outer iterations to obtain the final result.

%

\subsection{Reconstruction results}

We present in Table~\ref{tab:table2} the reconstruction results for the
dielectric constant of the scattering objects under consideration. The
second column of the table represents the dielectric constants at 3 GHz of
the scattering objects, independently measured by physicists from the
Department of Physics and Optical Science at the University of North
Carolina Charlotte (M. A. Fiddy and S. Kitchin). The percent number in the
brackets is the error in the measurement given by the standard deviation.
The measured relative permittivity of object 5 (a baseball) was not
available due to some technical issues.

Our reconstruction results and their relative errors compared with measured
dielectric constants are in the last two columns of the table. One
can observe a good accuracy for our reconstructions: only a few percent.
The dielectric constant of object 2 (wet wood) is relatively high which
leads to the worst error case in computations, as anticipated. We have
obtained best errors in objects 1 and 3. Note that object 1 is the
calibrating object. We also observe that the computational error in object 4
is 2.5 times less than its measurement error. The computational error does
not exceed a few percent in all four cases, just as in~\cite{Beili2012,
Thanh2014, Thanh2015}. \vspace{0cm} 
\begin{table}[h]
\caption{Measured and computed dielectric constant of the targets}
\label{tab:table2}\centering
\begin{tabular}{|c|c|c|c|}
\hline
Obj.\,ID & Measured $\varepsilon_r$\,(std.~dev.) & Computed $\varepsilon_r$
& Relative error \\ \hline
1 & 5.30 (1.6\%) & 5.44 & 2.6\% \\ \hline
2 & 8.48 (4.9\%) & 7.60 & 10.3\% \\ \hline
3 & 5.44 (1.1\%) & 5.55 & 2.0\% \\ \hline
4 & 3.80 (13.0\%) & 4.00 & 5.2\% \\ \hline
5 & not available & 4.76 & n/a \\ \hline
\end{tabular}%
\end{table}

We display in Figures~\ref{fi:6} and \ref{fi:7} the 3D-visualizations of
exact and reconstructed geometry of the first four targets using \texttt{%
isosurface} in Matlab. We do not display the image for object 5, which is
similar to that of object 4, since both are spherically shaped. In the
isosurface plotting the isovalue was chosen to be 50\% of the maximal value
of the computed $\varepsilon _{r}$. This choice is simply based on the
calibrating object and is applied to all other objects. From the figures we
can see that locations of the targets are well-reconstructed. It can also be
seen that although the shape reconstruction of targets is in general not
very good, the results for spherical shaped objects 3 and 4 seem to be
better than those of rectangular prism ones (objects 1 and 2), which is
reasonable.

%

\section{Summary}
\label{sect:summary}

We have presented in this paper the performance on a set of experimental
data of the recently developed globally convergent numerical method~\cite{Kliba2016} for a
Coefficient Inverse Problems for the Helmholtz equation.
This is a CIP with the backscatter multi-frequency data resulting from a
single direction of the incident plane wave, i.e., this is a CIP with the minimal amount of measurement. The key advantage of this method over traditional locally
convergent inversion techniques is that, under a
reasonable mathematical assumption in section~\ref{sect:initialTail}, it rigorously provides a good
approximation for the exact coefficient without \textit{a priori} knowledge of a
small neighborhood of this coefficient.

Our data were intentionally collected under non-ideal conditions: to model
the focus of our applications, which are detection and identification of
explosive devices. Thus, our data are contaminated by a significant amount
of noise. Denoising procedures are inapplicable here since our data have a
rich content. Thus, we have developed a new heuristic data preprocessing
procedure. The result of this procedure is used as the input for our
algorithm of section~\ref{sect:Algorithm1}. This procedure distills signals reflected by
targets of our interest from parasitic signals and also makes the
preprocessed data look similarly with the computationally simulated data.
In particular, one can see from a comparison of Figures~\ref{fi:3}(a) and~\ref{fi:3}(b) that the
data propagation step of this procedure indeed helps us separate the
signal reflected by our targets from those scattered by all
other objects in the room.

\begin{figure}[!h]
\centering
\subfloat[Exact location and shape of object
1]{\includegraphics[width=0.4\textwidth]{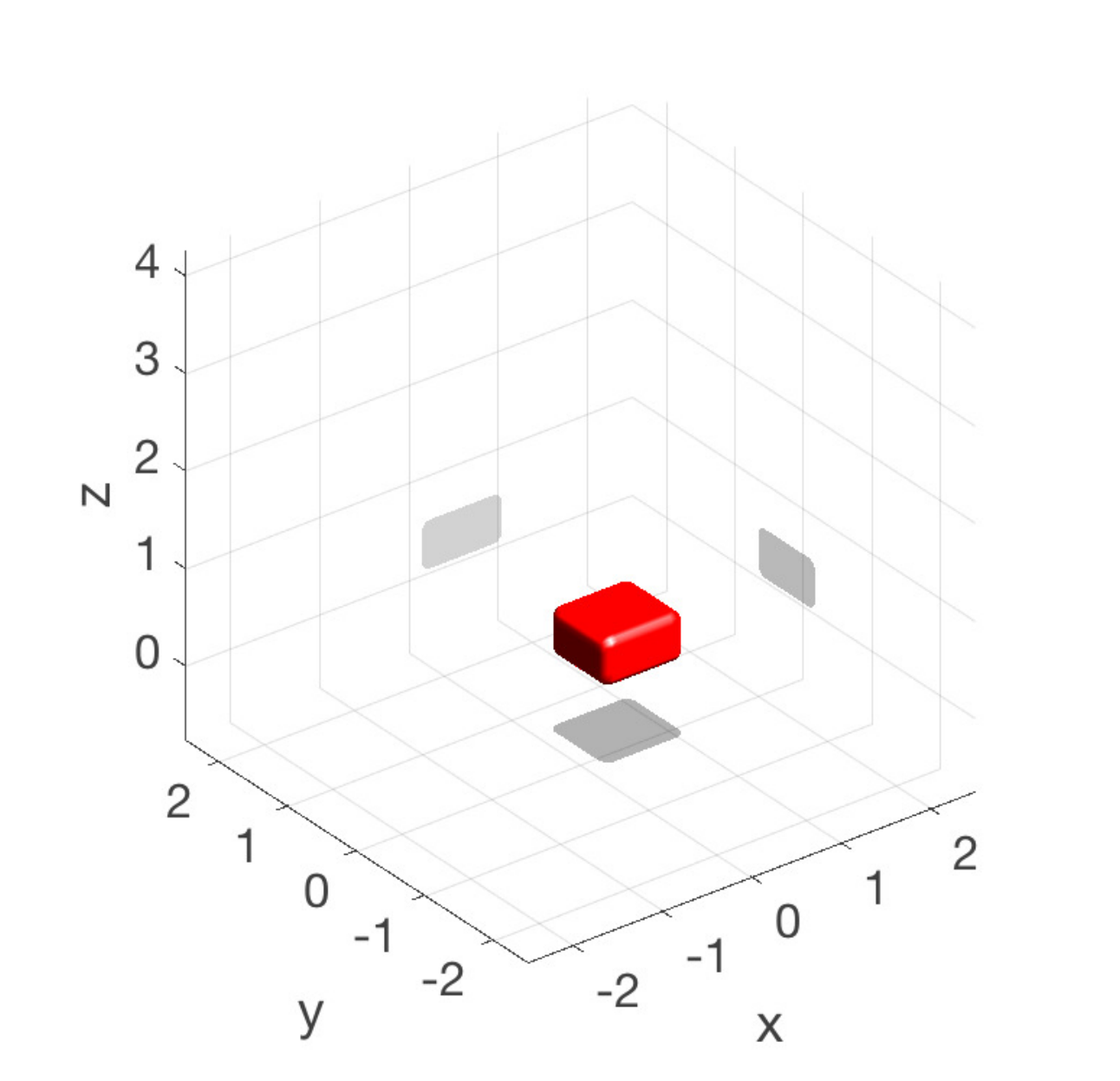}} \hspace{0.2cm} 
\subfloat[Reconstruction  result of object
1]{\includegraphics[width=0.4\textwidth]{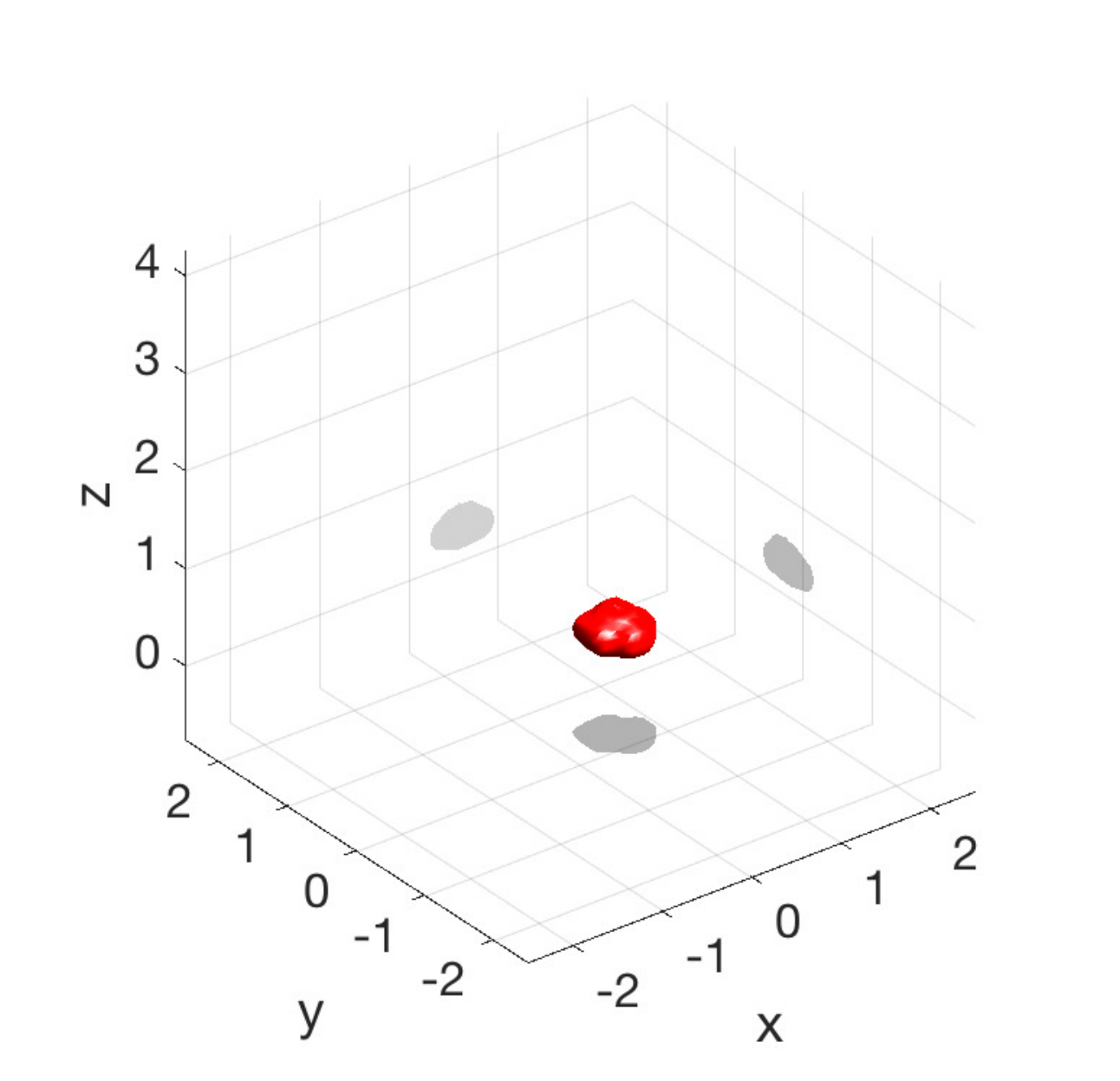}} 
\vspace{-0.0cm} 
\subfloat[Exact location and shape of object
2]{\includegraphics[width=0.4\textwidth]{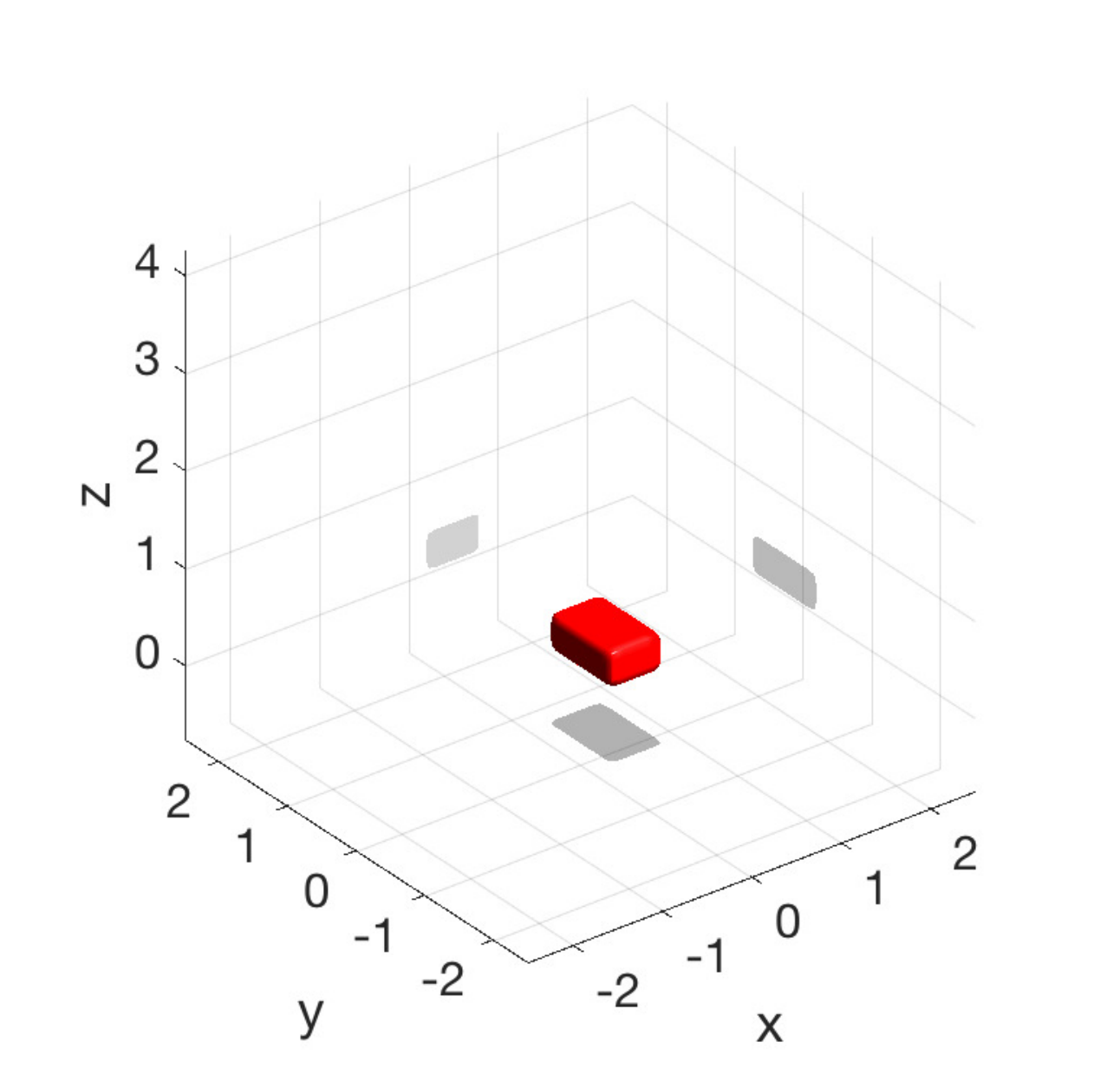}} \hspace{0.2cm} 
\subfloat[Reconstruction  result  of object
2]{\includegraphics[width=0.4\textwidth]{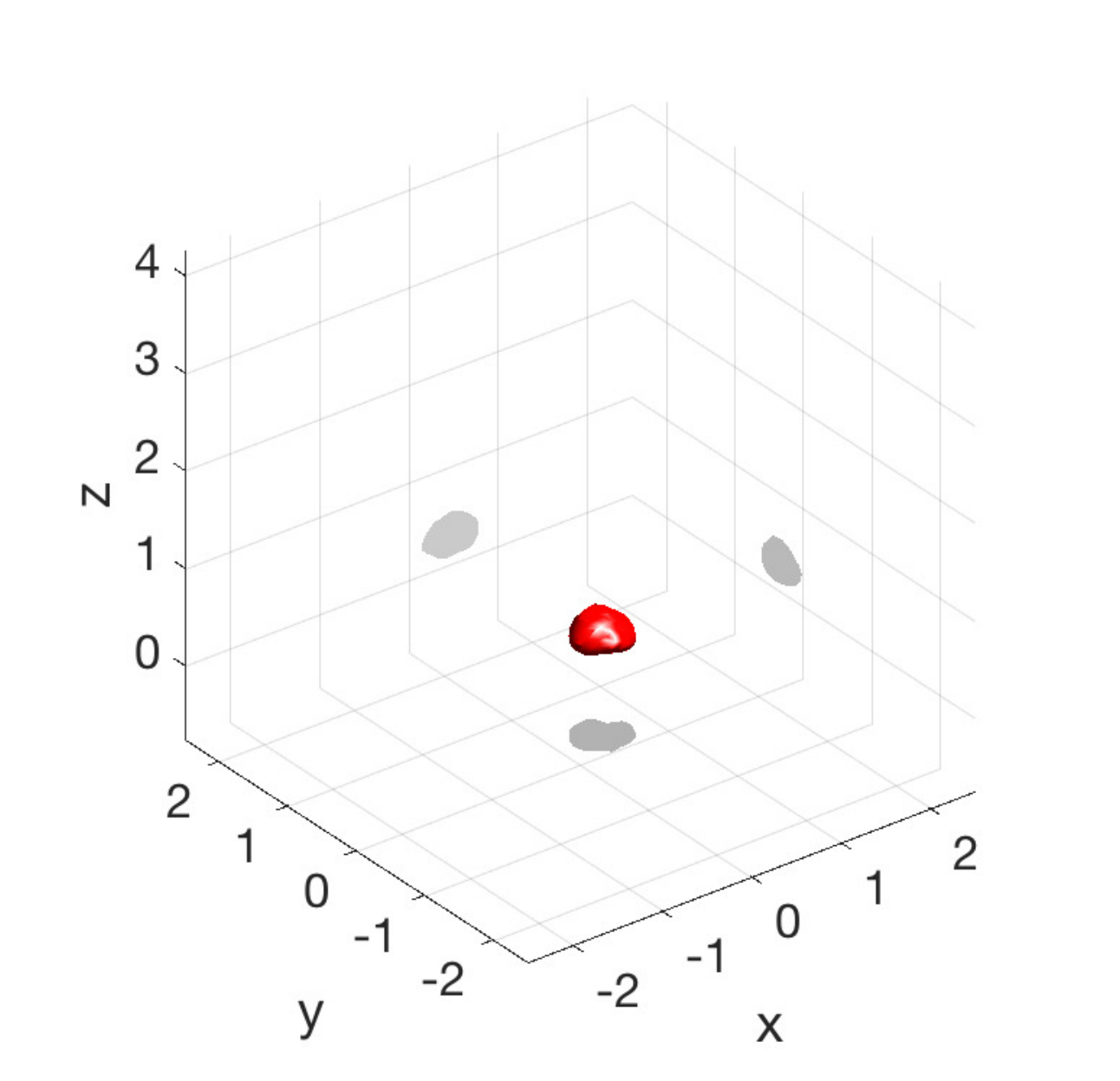}} 
\caption{ 3D visualizations of exact (left) and reconstructed (right)
geometry for objects 1 and 2.}
\label{fi:6}
\end{figure}

\begin{figure}[!h]
\centering
\subfloat[Exact location and shape of object
3]{\includegraphics[width=0.4\textwidth]{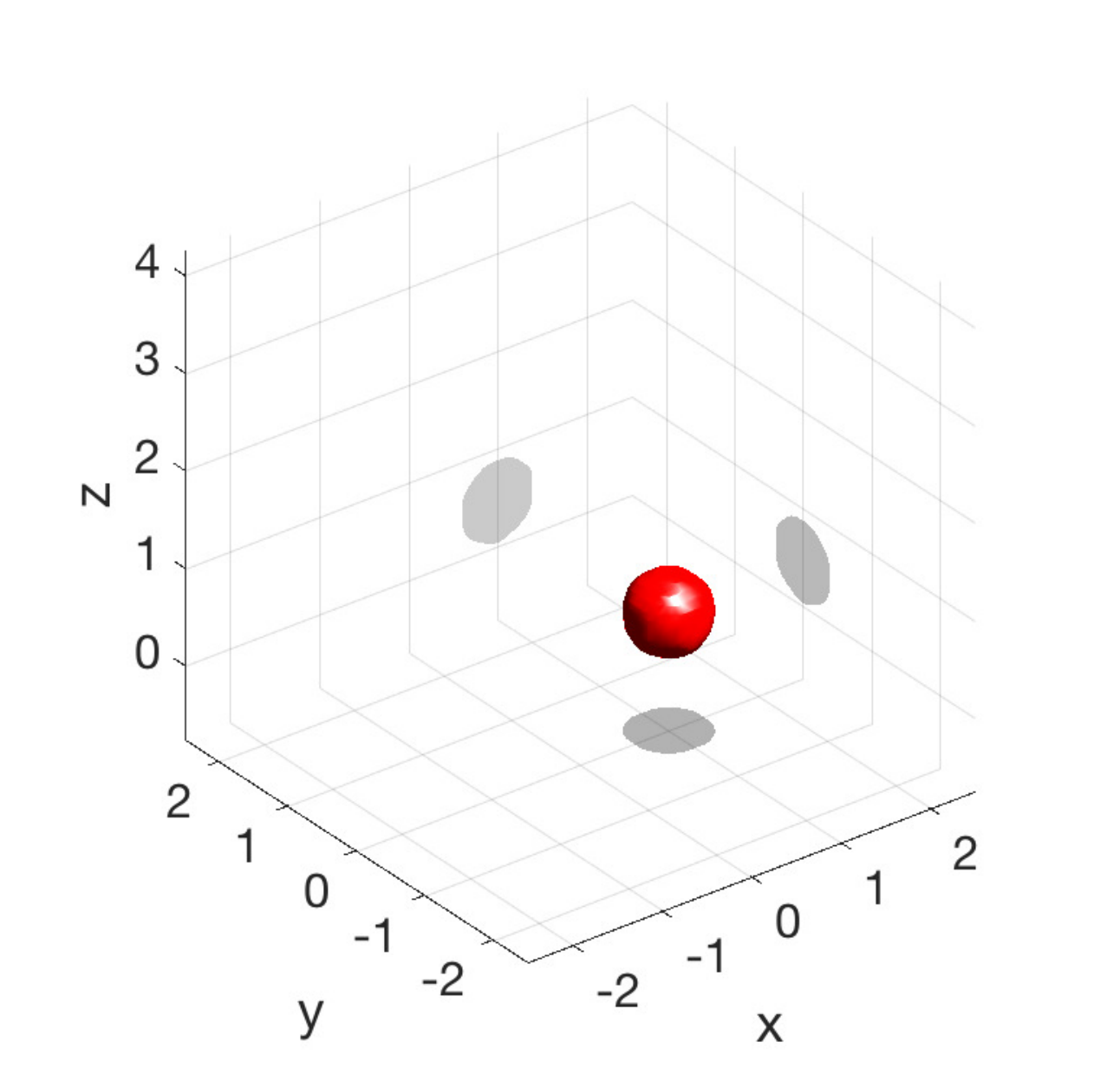}} \hspace{0.2cm} 
\subfloat[Reconstruction  result of object
3]{\includegraphics[width=0.4\textwidth]{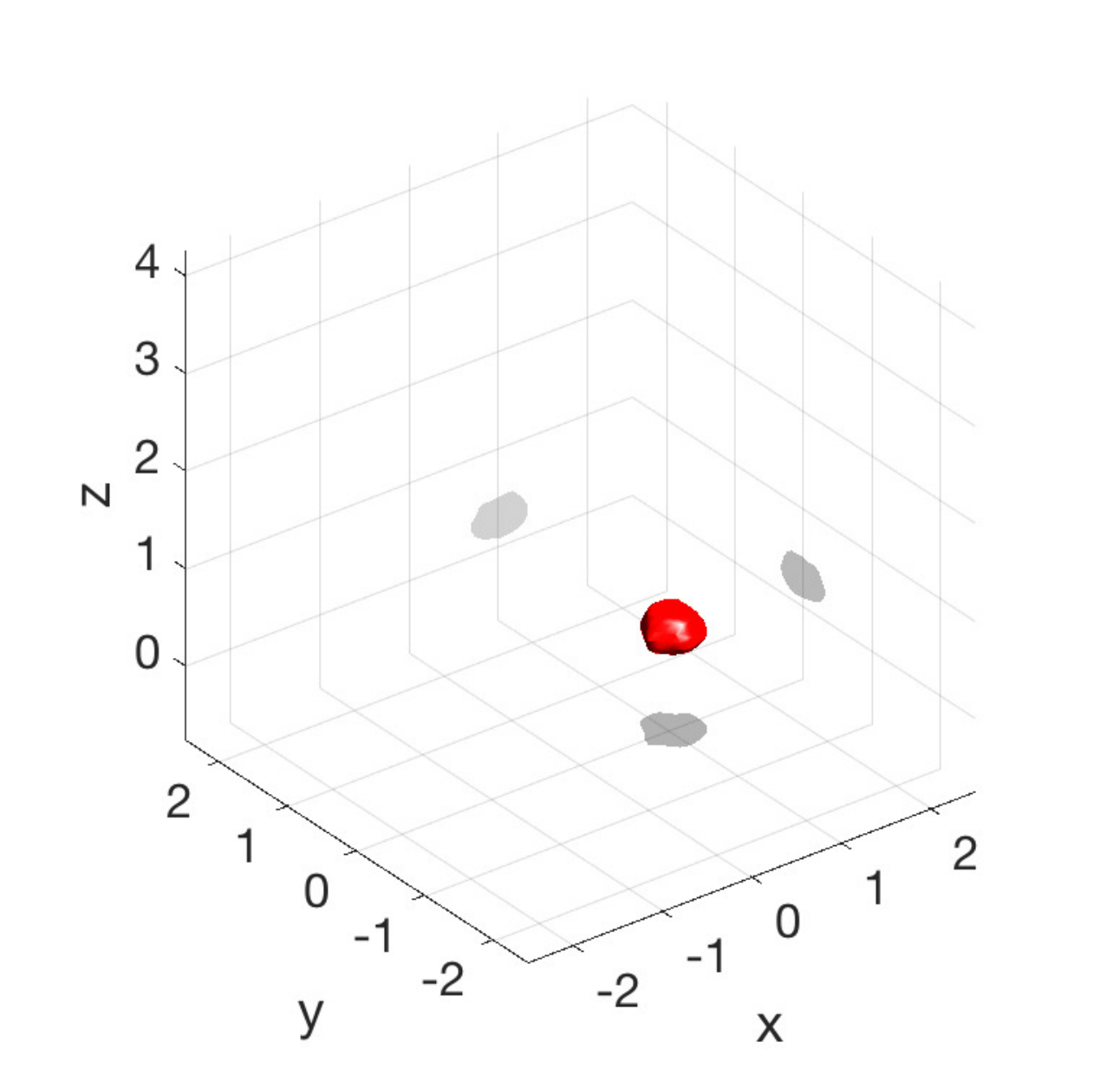}} 
\vspace{-0.0cm} 
\subfloat[Exact location and shape of object
4]{\includegraphics[width=0.4\textwidth]{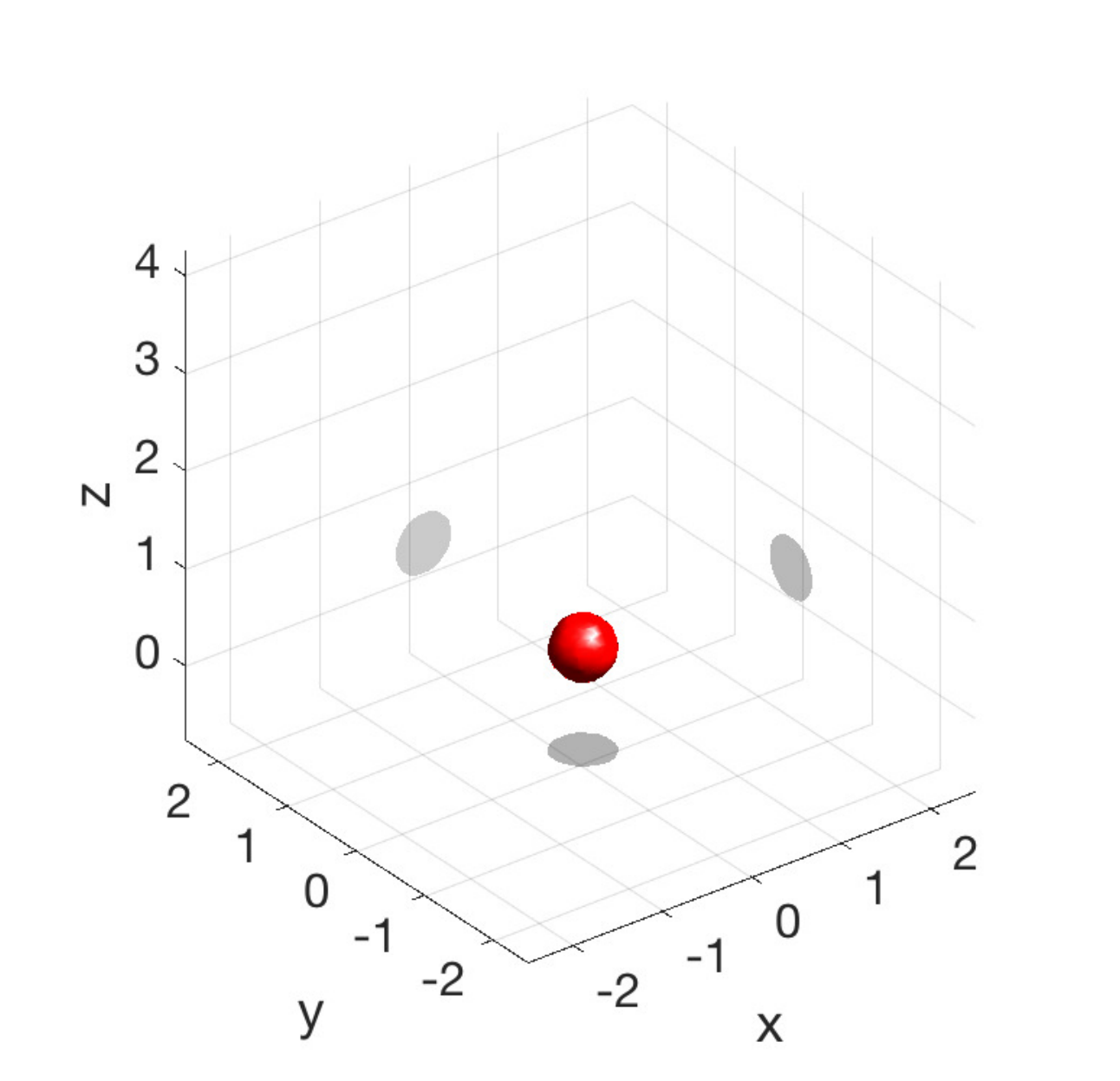}} \hspace{0.2cm} 
\subfloat[Reconstruction result of object
4]{\includegraphics[width=0.4\textwidth]{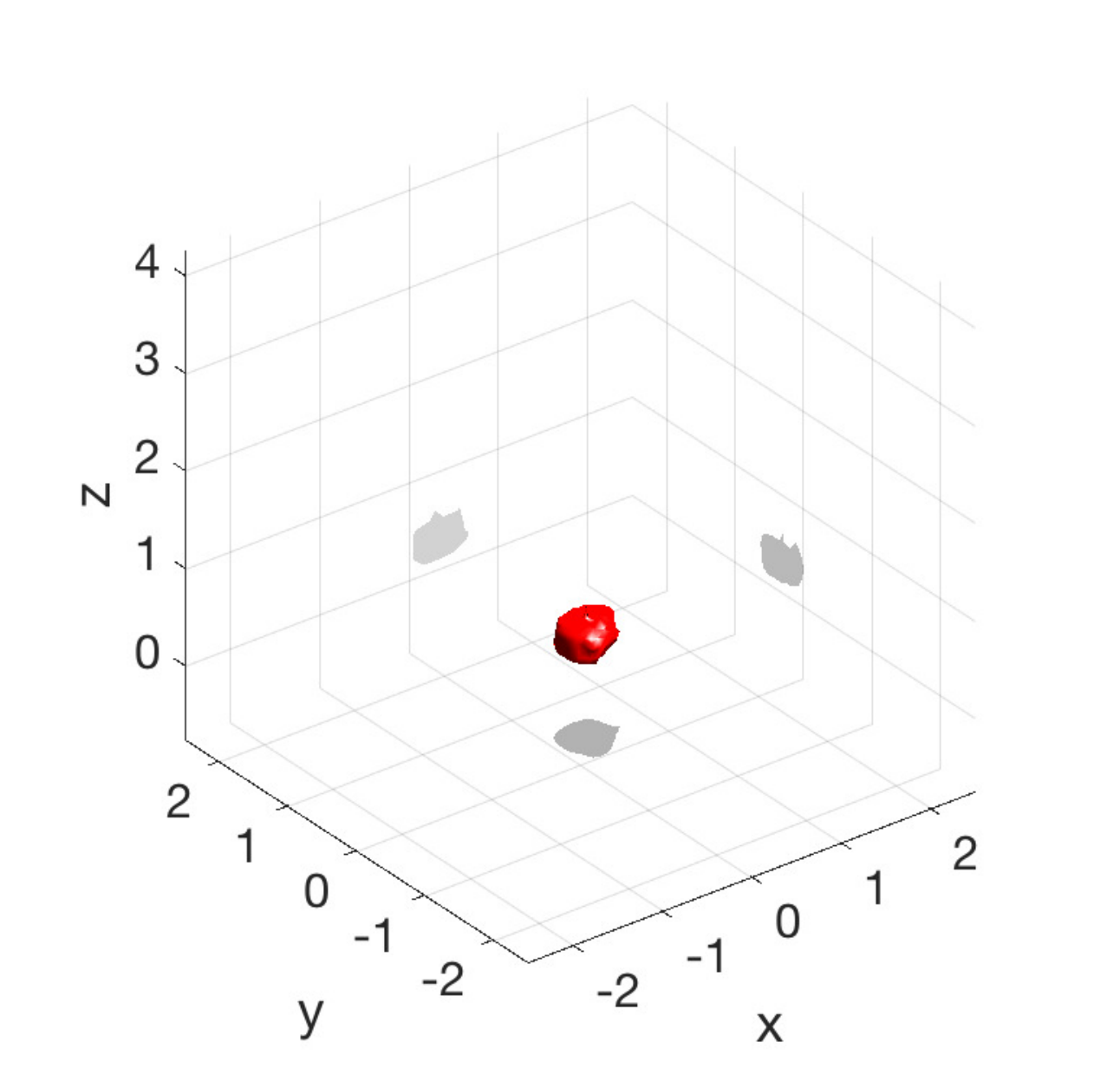}}
\caption{ 3D visualizations of exact (left) and reconstructed (right)
geometry for objects 3 and 4.}
\label{fi:7}
\end{figure}

Table \ref{tab:table2} shows that we come up with a good reconstruction accuracy
of dielectric constants of targets with the errors not exceeding a few
percent. One can also see from this table that our method can reconstruct
relatively large inclusion/background contrasts, the case that is well known
to be difficult for conventional least-squares approaches. Figures~\ref{fi:6} and~\ref{fi:7}
show that locations of targets are also reconstructed with a good accuracy.
This is achieved even though the data are highly noisy and we measured only
the backscatter data associated with a single incident field.

\section*{Acknowledgements}

This work was supported by US Army Research Laboratory and US Army Research
Office grant W911NF-15-1-0233 and by the Office of Naval Research grant
N00014-15-1-2330.
The authors are grateful to Mr. Steven Kitchen for his excellent work on
data collection.

%

\providecommand{\noopsort}[1]{}


\end{document}